\documentclass[11pt]{article}
\usepackage{amsmath,verbatim,hyperref,graphics}
\usepackage{epsfig,amssymb,latexsym}
\usepackage[all]{xy}

\numberwithin{equation}{section}
\newenvironment{proof}[1]{\addvspace{\medskipamount}\par\noindent{\it Proof#1}.}
{\unskip\nobreak\hfill$\Box$\par\addvspace{\medskipamount}}

\def\q{\quad}

\def\td{\mathrm{d}}
\def\BEN{\begin{enumerate}}  \def\BI{\begin{itemize}}
\def\EEN{\end{enumerate}}   \def\EI{\end{itemize}}

\def\beq{\begin{eqnarray}} \def\eeq{\end{eqnarray}}

\def\al*#1{\begin{align*}#1\end{align*}}

\def\ga*#1{\begin{gather*}#1\end{gather*}}

\def\alat*#1#2{\begin{alignat*}{#1}#2\end{alignat*}}
\def\bea{\begin{eqnarray*}}
\def\eea{\end{eqnarray*}}
\def\ml*#1{\begin{multline*}#1\end{multline*}}

\def\mbb{\mathbb} \def\mbf{\mathbf} 
\def\mc{\mathcal}  \def\ovl{\overline}

\def\F{\mathcal{F}}
\def\P{{\mathbb P}} \def\le{\left} \def\ri{\right} 
\def\E{{\mathbb E}}  \def\Z{{\mathbb Z}} \def\R{{\mathbb R}}
 \def\Q{{\mathbb Q}}

\def\G{\Gamma}

\def\te#1{\mathrm{e}^{#1}}   

\def\WH{\widehat}

  \def\im{\item} 
 \def\U{\mathbb U}    
     
\def\e{\varepsilon}

\def\w{\omega} \def\q{\qquad} \def\D{\Delta}
   
  \def\td{\text{\rm d}}
 
\numberwithin{equation}{section}

\numberwithin{equation}{section}
\newtheorem{Thm}{Theorem}[section]
\newtheorem{Cor}[Thm]{Corollary}
\newtheorem{Def}[Thm]{Definition}
\newtheorem{Prop}[Thm]{Proposition}
\newtheorem{Lemma}[Thm]{Lemma}
\newtheorem{remarks}[Thm]{Remark}
\newtheorem{Cond}[Thm]{Condition}

\newtheorem{as}[Thm]{Assumption}

\newtheorem{examples}[Thm]{Example}
\newtheorem{foo}[Thm]{Remarks}
\newenvironment{Example}{\begin{examples}\rm}{\end{examples}}

\newenvironment{Remark}{\begin{remarks}\rm}{\end{remarks}}
\newenvironment{Rem}{\begin{remarks}\rm}{\end{remarks}}

\newenvironment{As}{\begin{as}\rm}{\end{as}}

\begin{document}

\title{{\bf\Large On dynamic spectral risk measures, a limit theorem\\ and optimal portfolio allocation}}
\author{Dilip Madan\footnote{Robert H. Smith School of Business,
University of Maryland, dbm@rhsmith.umd.edu}\qquad 
Martijn Pistorius\footnote{(Corr. author) Department of Mathematics, Imperial College London, m.pistorius@imperial.ac.uk}
\qquad
Mitja Stadje%
\footnote{Faculty of Mathematics and Economics,
Universit\"{a}t Ulm, mitja.stadje@uni-ulm.de
\newline
{\em Keywords:} Spectral risk measure, dynamic risk measure, $g$-expectation, Choquet expectation, distortion, 
(strong) time-consistency, limit theorem, dynamic portfolio optimisation, (extendend) HJB equation, 
sub-game perfect Nash equilibrium, dynamic Markowitz-type portfolio problem.
\newline
{\em Subject Classification.} {Primary: 60H10, 60Fxx; Secondary: 60H99, 91B30.
\newline
{\em JEL classification}: G32
}
}
}
\date{}

\maketitle
\vspace{-0.6cm}
%\begin{spacing}{0.83}
\begin{quote}
{\small{\bf Abstract.} In this paper we propose the notion of continuous-time
{\em dynamic spectral risk-measure} (DSR). Adopting a Poisson random measure setting, 
we define this class of dynamic coherent risk-measures in terms of 
certain backward stochastic differential equations. By establishing a functional limit theorem,
we show that DSRs may be considered to be (strongly) time-consistent
continuous-time extensions of iterated spectral risk-measures, which are obtained 
by iterating a given spectral risk-measure (such as Expected Shortfall) along a 
given time-grid. Specifically, we demonstrate that any DSR arises in the limit of a sequence 
of such iterated spectral risk-measures driven by lattice-random walks, under suitable scaling 
and vanishing time- and spatial-mesh sizes. 
To illustrate its use in financial optimisation problems, 
we analyse a dynamic portfolio optimisation problem under a DSR. 
}
\end{quote}
%\end{spacing}
%\vspace{-0.3cm}

\section{Introduction}
\label{intro}
Financial analysis and decision making rely on quantification and modelling of future risk exposures. 
A systematic approach for the latter was put forward in \cite{ADEH}, laying the foundations
of an axiomatic framework for coherent measurement of risk. 
A subsequent breakthrough was the development and application of the notion of 
backward stochastic differential equations (BSDEs) 
in the context of risk analysis, which gave rise to the (strongly) time-consistent extension 
of coherent risk-measures to continuous-time dynamic settings \cite{PengR1,RG}. 
Building on these advances, we consider in this article
a new class of such continuous-time dynamic coherent risk measures,
which we propose to call {\em dynamic spectral risk measures} (DSRs). 

Quantile-based coherent risk measures, such as {\em Expected Shortfall}, belong to the most widely used risk-measures 
in risk analysis, and are also known as {\em spectral risk measures}, {\em Choquet expectations} (based on probability distortions) and {\em Weighted VaR}; see \cite{Ac,CWV,Kusuoka,Wang96}.
In order to carry out, for instance, an analysis of portfolios involving dynamic rebalancing, 
one is lead to consider the (strongly) time-consistent extension of such coherent risk-measures 
to given time-grids, which are defined by iterative application of the spectral risk-measure along these particular 
grids. Due to its continuous-time domain of definition a DSR  
is, in contrast, independent of a grid structure. While the latter holds for any continuous-time risk measure
we show that DSRs emerge as the limits of such iterated spectral risk measures when the time-step vanishes and under appropriate scaling of the parameters, by establishing a functional limit theorem. 

To explore its use in financial decision problems, we consider subsequently a dynamic portfolio optimisation problem under DSR, 
which we analyse in terms of its associated  Hamilton-Jacobi-Bellman (HJB) equation. In the case of a 
long-only investor (who is allowed neither to borrow nor to short-sell stocks) we identify explicitly dynamic optimal allocation strategies.

% Time-consistency 

DSR, like any dynamic risk-measure obtained from a BSDE, is {\em (strongly) time-consistent} 
in the sense that if the value of a random variable $X$ is not larger than $Y$ under DSR at time $t$ almost surely,
then the same relation holds at earlier times $s$, $s<t$. For dynamic risk-measures 
the property of strong time-consistency is well known to be equivalent to {\em recursiveness}, a tower-type property which is referred to as filtration-consistency in \cite{CHMP}  
Such concepts have been investigated extensively in the literature; among others we mention 
~\cite{ADEHK,CDK,cohenelliota,D2006,FollmerP,JR,KlSchw,Riedel}. For studies on weaker forms of time-consistency we 
refer to \cite{RoorS,Tutsch,Weber}. 
%
% Time-consistency2
%

The notion of strong time-consistency in {\em economics} goes back at least as far
as \cite{Strotz}  and has been standard in the economics literature ever since; see for instance 
\cite{CE,DE1992,ES2003,EZ1989,HS2008,K1960,KP1978}. 

Due to their recursive structure financial optimisation problems, such as utility opimisation under the
entropic risk-measure and related robust portfolio optimisation problems satisfy the Dynamic Programming Principle  and admit time-consistent dynamically optimal strategies (see for instance \cite{B,R} and references therein). 
In Section~\ref{sec:DynOpt} we demonstrate that this also holds for the optimal portfolio allocation problem 
phrased in terms of the minimisation under a DSR, and phrase and solve this problem via the associated HJB equation.   

% Limit theorem + Static spectral risk-measures

For a given DSR, the functional limit 
theorem that we obtain (see Theorem~\ref{main}) shows how to construct 
an approximating sequence of iterated spectral risk-measures 
driven by lattice random walks, suggesting an effective method to evaluate functionals under a given DSR 
and solutions to associated PIDEs, by recursively applying (distorted) Choquet expectations. 
The functional limit theorem involves a certain non-standard scaling 
of the parameters of the iterated spectral risk measures, which is given in Definition~\ref{def:PsiD}.
The advantage of this approximation method is that it sidesteps 
the (typically non-trivial) task of computing the Malliavin derivatives. 
A numerical study is beyond the scope of the current paper, and is left for future research. 

While one may prove the functional limit theorem directly 
through duality arguments, we present in the interest of brevity a proof that draws on the convergence results obtained 
in \cite{MPS1} for weak approximation of BSDEs. In the literature various related 
convergence results are available, of which we next mention a number (refer to \cite{MPS1} for additional references). The construction of continuous-time dynamic risk-measures arising as limits of discrete-time ones was studied 
in \cite{SIME} in a Brownian setting. In a more general setting including in addition finitely many Poisson processes,
\cite{KuMo} presents a limit theorem for recursive coherent quantile-based risk measures, which is proved via an associated non-linear partial differential equation. In \cite{DNM} a Donsker-type theorem is established under a $G$-expectation.  
\bigskip

{\bf Contents.} The remainder of the paper is organised as follows. 
In Section~\ref{sec:gexpF} we collect preliminary results 
concerning dynamic coherent risk measures and related BSDEs, adopting a pure jump setting driven by a Poisson random measure. 
In Section~\ref{sec:ChE} we are concerned with the Choquet-type integrals which appear in the definitions of 
dynamic and iterated spectral risk measures. With these results in hand, we phrase 
the definition of a DSR in Section~\ref{sec:DSR} and identify its dual representation. In Section~\ref{sec:lim} we present the functional limit theorem for iterated spectral risk measures. Finally, in Section~\ref{sec:DynOpt} we turn to the study of a dynamic portfolio allocation problems under a DSR.

\section{Preliminaries}\label{sec:gexpF}
In this section we collect elements of the theory of  
time-consistent dynamic coherent risk measures and associated BSDEs, in both continuous-time
and discrete-time settings. To avoid repetition we state some results and definitions 
in terms of the index set $\mathcal I$, which is taken to be either $\mathcal I=[0,T]$ or
$$\mathcal I = \pi_\Delta := \{t_i=i\Delta, i=0,\ldots, N\},\  \text{with}\ \Delta = T/N$$
for some $N\in\mathbb N$ and  $T>0$.

\subsection{Time-consistent dynamic coherent risk measures}
On some filtered probability space $(\Omega,\mc F, \mbf F, \P)$ with $\mbf F=(\mc F_t)_{t\in I}$, 
we consider risks described by random variables $X\in\mc L^p = \mc L^p(\mc F_T)$, 
$p>0$, the set of  $\mc F_T$-measurable random variables $X$ with $\E[|X|^p] = \int_\Omega |X|^p\td P<\infty$. 
We denote by $\mc L^p_t = \mc L^p(\mc F_t)$ and $\mc L^p(\mc G)$ the elements 
$X$ in $\mc L^p(\mc F)$ that are measurable with respect to the sigma-algebras $\mc F_t$ and  $\mc G\subset\mc F$, respectively,
and by $\mc L^\infty, \mc L^\infty_t, \mc L^\infty(\mc G)$ the collections of bounded elements in $\mc L^p$, 
$\mc L^p_t$ and $\mc L^p(\mc G)$. Let $\mc S^2(\mathcal I)$ denote the space of $\mbf F$-adapted semi-martingales $Y = (Y_t)_{t\in\mathcal I}$ that are square-integrable in the sense that $\|Y\|^2_{\mc S^2(\mathcal I)}<\infty$, where 
$$\|Y\|^2_{\mc S^2(\mathcal I)}:=\E\le[\sup_{t\in\mathcal I}|Y_t|^2\ri].$$ 
For a given measure $\mu$ on a measurable space 
$(\mathbb U,\mc U)$ we denote by $\mc L^p(\mu)$, $p>0$, 
the set of Borel functions $v: \mbb U\to\mathbb R$ 
with $|v|_{p,\mu}<\infty$, where
$$|v|_{p,\mu} := \left(\int_{\mbb U} |v(x)|^p \mu(\td x)\right)^{1/p},
$$
and by $\mc L_+^p(\mu)$ the set of non-negative elements in $\mc L^p(\mu)$.

Dynamic coherent risk-measures and (strong) time-consistency, we recall, are defined as follows in an $\mc L^2$-setting:

\begin{Def}\label{def:rho}A {\em  dynamic coherent 
		risk measure} $\rho=(\rho_t)_{t\in I}$ is a map $\rho:\mc L^2\to  \mc S^2(\mathcal I)$
	that satisfies the following properties:
	\BI
	\im[{\bf (i)}] {\em (cash invariance)} for $m\in \mc L^2_t$, $\rho_t(X+m) = \rho_t(X) - m$; 
	\im[{\bf (ii)}] {\em (monotonicity)} for $X,Y\in \mc L^2$ with $X\ge Y$, $\rho_t(X)\leq\rho_t(Y)$; 
	\im[{\bf (iii)}] {\em (positive homogeneity)} for $X\in \mc L^2$ and $\lambda\in \mc L^\infty_t,$ 
	$\rho_t(|\lambda| X) = |\lambda|\rho_t(X)$;
	\im[{\bf (iv)}] {\em (subadditivity)} for $X,Y\in \mc L^2$,
	$\rho_t(X + Y) \leq \rho_t(X) + \rho_t(Y)$. 
	\EI
\end{Def} 
\begin{Def}\label{def:coh}
	A dynamic coherent risk measure $\rho$ is called  {\em (strongly) time-consistent} if either of the following holds: 
	\BI
	\im[{\bf (v)}] {\em (strong time-consistency)} for $X,Y\in \mc L^2$ and $s,t$ with $s\leq t$,\ 
	$\rho_t(X) \leq \rho_t(Y) \Rightarrow \rho_s(X) \leq \rho_s(Y)$; 
	\im[{\bf (vi)}] {\em (recursiveness)} for $X\in \mc L^2$ and $s,t$ with $s\leq t$, $\rho_s(\rho_t(X)) = \rho_s(X)$. 
	\EI
\end{Def}

\noindent For a proof of the equivalence of items {\bf (v)} and {\bf (vi)} we refer to F\"{o}llmer and Schied (2011); 
for a discussion of (the unconditional version of) the properties {\bf (i)}--{\bf (iv)}
see \cite{ADEH,ADEHK}. 
One way to construct a time-consistent  dynamic coherent risk measure is as solution to an associated backward stochastic differential equation (BSDE) or backward stochastic difference equation (BS$\Delta$E). 
To ensure that such dynamic risk measures satisfy properties {\bf (i)}--{\bf (iv)}  
the corresponding driver functions are to be positively homogeneous, subadditive and should not dependent on the value of the risk-measure (see Proposition~11 in \cite{RG} and Lemma 2.1 in \cite{CHMP}.
For background on the notion of strong time-consistency and its relation to $g$-expectations we refer to 
\cite{BN2,BN1,PengR1,RG}. Specifically,  in our setting 
such driver functions are defined as follows:

\begin{Def}\label{def:cohg} For a given Borel measure $\mu$ on $\mathbb R^k\backslash\{0\}$ we call a function 
	$g:\mathcal I \times \mc L^2(\mu)\to\mbb R$ a {\em driver function}
	if for any $z\in\mc L^2(\mu)$ $t\mapsto g(t, z)$ is continuous (in case $\mathcal I = [0,T]$) 
	and the following holds:
	\BI
	\im[{\bf (i)}] {\em (Lipschitz-continuity)} for some $K\in\mbb R_+\backslash\{0\}$ and 
	any $t\in\mathcal I$ and $ z_1, z_2\in \mc L^2(\mu)$
	\begin{equation*}%\label{Klip}
	|g({t}, z_1) - g({t}, z_2)| \leq K | z_1- z_2|_{2,\mu}.
	\end{equation*}
	\EI
	A driver function $g$ is called {\em coherent} if the following hold:
	\BI
	\im[{\bf (ii)}] {\em (positive homogeneity)} for any $r\in\mbb R_+$, $t\in\mathcal I$ 
	and $ z\in \mc L^2(\mu)$, we have  
	$$g({t}, r  z) =  
	r g({t},  z);$$
	\im[{\bf (iii)}] {\em (subadditivity)} for any $t\in\mathcal I$
	and $ z_1, z_2\in \mc L^2(\mu)$, we have  
	$$g({t},   z_1 +  z_2) \leq  g({t}, z_1) + 
	g({t}, z_2).$$
	\EI
\end{Def}
We describe next the dynamic coherent risk-measure defined via 
the  BSDEs (if $\mathcal I=[0,T]$) or BS$\Delta$Es 
(if $\mathcal I$ is a finite partition of $[0,T]$) corresponding to coherent 
driver function functions. 

\subsection{Discrete-time lattice setting}
We turn first to the discrete-time lattice setting, fixing a uniform partition $\pi=\pi_\Delta$ of $[0,T]$ 
with as before $\Delta = T/N$ for some $N\in\mathbb N$. Let $L^{(\pi)}=(L^{(\pi)}_t)_{t\in\pi}$ denote 
a square-integrable zero-mean random walk starting at zero and taking values in 
$(\sqrt{\Delta}\mathbb Z)^k$,  and let $\mathbf F^{(\pi)} = (\mathcal F^{(\pi)}_t)_{t\in\pi}$ denote 
the filtration generated by $L^{(\pi)}$. 
Furthermore, we let $g^{(\pi)}$ be a coherent driver function as in Definition~\ref{def:cohg} 
with $\mathcal I=\pi$ and $\mu(\td x)$ equal to the scaled law $\nu^{(\pi)}(\td x)$ of 
$\Delta L^{(\pi)}_t = L^{(\pi)}_{t + \Delta} - L^{(\pi)}_{t}$, $t\in\pi\backslash\{T\}$, 
given by
\begin{equation}\label{eq:nupi}
\nu^{(\pi)}(\td x) := \frac{1}{\Delta}\,\P(\Delta L_{t}^{(\pi)}\in\td x),\q x\in (\sqrt{\Delta}\mathbb Z)^k.
\end{equation}
Since the predictable representation property continues to hold in this setting, 
the BS$\Delta$E for $(Y^{(\pi)}, Z^{(\pi)})$ 
corresponding to final value $- X^{(\pi)}\in \mc L^2(\F_T^{(\pi)})$ and driver function $g^{(\pi)}$
takes the following form, which is analogous to the one in 
continuous-time case given in \eqref{bsde1} below:
\begin{eqnarray}\nonumber
Y_t^{(\pi)} &=& - X^{(\pi)} + \sum_{s=t}^{T-\Delta} g^{(\pi)}(s,  Z_s^{(\pi)})\Delta \\
&&- \sum_{s=t}^{T-\Delta} \left( Z^{(\pi)}_s(\Delta L^{(\pi)}_s)I_{\{\Delta L^{(\pi)}_s\neq 0\}} - 
\E\left[\left. Z^{(\pi)}_s(\Delta L^{(\pi)}_s)I_{\{\Delta L^{(\pi)}_s\neq 0\}}\right| \mathcal F^{(\pi)}_s
\right]\right)
\label{BSDD}
\end{eqnarray}
for $t\in\pi\backslash\{T\}$ and with $Y_T^{(\pi)} = - X^{(\pi)}$,
where $I_A$ denotes the indicator of a set $A$.
In difference notation the BS$\Delta$E \eqref{BSDD} is for $t\in\pi\backslash\{T\}$ given by
\begin{eqnarray}\nonumber
\!\!\!\!\!\!\!\Delta Y_t^{(\pi)} &=& -g^{(\pi)}(t,  Z_t^{(\pi)})\Delta  \\
&+&  Z^{(\pi)}_t(\Delta L^{(\pi)}_t)I_{\{\Delta L^{(\pi)}_t\neq 0\}} - 
\E\left[\left. Z^{(\pi)}_t(\Delta L^{(\pi)}_t)I_{\{\Delta L^{(\pi)}_t\neq 0\}}\right| \mathcal F^{(\pi)}_t
\right]
\label{BSDD1}
\end{eqnarray}
with $Y_T^{(\pi)} = - X^{(\pi)}$.
A pair $(Y^{(\pi)},  Z^{(\pi)})$ is a solution of the BS$\Delta$E 
if,  for any $t\in\pi$, it satisfies \eqref{BSDD} with 
$$Y_t^{(\pi)}\in \mathcal L^2(\F^{(\pi)}_t),\q 
Z^{(\pi)}_t \in 
\mathbb L^2_t:= \mathcal L^2(\nu^{(\pi)}(\td x)\times\td\P,\mc B((\sqrt{\Delta}\mbb Z)^k)\otimes \F_t^{(\pi)}).
$$
If the Lipschitz-constant $K=K^{(\pi)}$ of the driver function $g^{(\pi)}$ is strictly smaller than the reciprocal $1/\Delta$ of the mesh-size then it follows from \cite[Propositions 3.1 and 3.2]{MPS1} 
that there exists a unique solution 
$(Y^{(\pi)},  Z^{(\pi)})$
to the BS$\Delta$E  
which satisfies the following relations for $t\in\pi$:
\begin{eqnarray}
\label{Ypi}
Y_t^{(\pi)} &=& g^{(\pi)}(t, Z_t^{(\pi)})\Delta + \E\left[\left.Y^{(\pi)}_{t+\Delta}\right|\F_t^{(\pi)}\right],\\ 
\nonumber
Z_t^{(\pi)}(x) &=& \E\left[\left.Y^{(\pi)}_{t+\Delta}\right|\F_t^{(\pi)}\vee\{\Delta L^{(\pi)}_t = x\}\right]\\ \label{tildeZpi}
&&- \E\left[\left.Y^{(\pi)}_{t+\Delta}\right|\F_t^{(\pi)}\vee\{\Delta L^{(\pi)}_t = 0\}\right]
\end{eqnarray}
for $x\in(\sqrt{\Delta}\mathbb Z)^k$,
where $\F_t^{(\pi)}\vee\{\Delta L^{(\pi)}_t = x\} := \F_t^{(\pi)}\vee\sigma(\{\Delta L^{(\pi)}_t = x\})$ 
denotes the smallest sigma-algebra containing $\F_t^{(\pi)}$ as well as 
the sigma-algebra $\sigma(\{\Delta L^{(\pi)}_t = x\})$ generated by $\{\Delta L^{(\pi)}_t = x\}$.   
In analogy with the continuous-time case (reviewed below), 
the dynamic coherent risk-measure associated to the solution to the BS$\Delta$E 
is defined as follows:
\begin{Def}\label{def:disrim}
	For a coherent driver function $g^{(\pi)}$ as in Definition~\ref{def:cohg} 
	with $\mathcal I=\pi$ and $\mu(\td x)=\nu^{(\pi)}$ and the solution 
	$(Y^{(\pi)}, Z^{(\pi)})$ of the corresponding BS$\Delta$E \eqref{BSDD}, 
	$\rho^{g^{(\pi)},(\pi)}=(\rho^{g^{(\pi)},(\pi)}_{t})_{t\in\pi}$  
	denotes the dynamic coherent risk measure given by $\rho^{g^{(\pi)},(\pi)}_{t}: \mc L^2(\mc F^{(\pi)}_T)\to \mc L^2(\mc F^{(\pi)}_t)$
	with
	$$\rho^{g^{(\pi)},(\pi)}_t(X) = Y_t^{(\pi)}.$$ 
\end{Def}

\subsection{Continuous-time setting}\label{sec:cont}
In the continuous-time case ($\mathcal I=[0,T]$) we consider risky positions described by random variables $X$ that 
are measurable with respect to $\mathcal F_T$, where $\mbf F=\{\mathcal F_t\}_{t\in[0,T]}$ denotes the right-continuous and 
completed filtration generated by a Poisson random measure
$N$ on $[0,T]\times\mbb R^{k}\backslash\{0\}$ for some $k\in\mathbb N$. 
We suppose throughout that the associated L\'{e}vy measure $\nu$ satisfies the following condition:
\begin{As}\label{As:LevyC}
	The L\'{e}vy measure $\nu$ associated to the Poisson random measure 
	$N$ has no atoms and, for some $\varepsilon_0>0$, 
	$\nu_{2+\varepsilon_0}\in \mbb R_+\backslash\{0\}$ where for $p\ge 0$
	$$\nu_{p}:=\int_{\mbb R^{k}\backslash\{0\}}|x|^{p}\nu(\td x) .$$
\end{As}
We denote by $\tilde N(\td t\times\td x)= N(\td t\times\td x) - \nu(\td x)\td t$ the compensated Poisson random measure 
and by $L=(L_t)_{t\in[0,T]}$ the (column-vector) L\'{e}vy process given by 
$$L_t=\int_{[0,t]\times\mathbb R^k\backslash\{0\}}x \tilde N(\td s\times\td x).$$ 
Under Assumption~\ref{As:LevyC} we have $\E[|L_t|^{2+\varepsilon_0}]<\infty$ for any $t\in[0,T]$ (see \cite[Theorem 25.3]{Sato}

Let $\tilde{\mc H}^2$ denote the set of $\tilde{\mc P}$-measurable square-integrable processes,
where, with $\mc P$ denoting the predictable sigma-algebra, $\tilde{\mc P}=\mc P\otimes\mc B(\mbb R^k\backslash\{0\})$, 
and let $\mathcal U$ denote the Borel sigma-algebra induced by the $L^2(\nu(\td x))$-norm.
In particular, $U\in \tilde{\mc H}^2$ is such that  
$\|U\|_{\tilde{\mc H}^2} <\infty$, where 
$$
\|U\|_{\tilde{\mc H}^2} := \E\left[\int_0^T|U_t|^2_{2,\nu}\td t\right].
$$
Moreover, let $\mathcal M_2$ denote the set of probability measures 
$\mathbb Q = \mathbb Q^\xi$ on $(\Omega,\F_T)$ that are absolutely continuous with respect to $\P$ 
with square-integrable Radon-Nikodym derivatives $\xi\in \mathcal L^2_+(\F_T)$, and write $\mathcal S^2:= \mathcal S^2[0,T]$. 
%\noindent

Let us next consider a coherent driver function $g$ as in Definition~\ref{def:cohg} with $\mu=\nu$ and $\mathcal I=[0,T]$.
and fix a final condition $X\in\mathcal L^2$.
The associated BSDE for the pair $(Y,Z)\in\mc S^2\times\tilde{\mc H}^2$ is 
given by 
\begin{eqnarray}\label{bsde1}
Y_t = -X + \int_t^T g(s, Z_s)\,\td s  - \int_{(t,T]\times\mbb R^{k}\backslash\{0\}}Z_s(x)\,\tilde N(\td s\times\td x)\end{eqnarray}
for $t\in[0,T]$. This BSDE, 
we recall from \cite{BBP}, admits a unique solution.
By combining \cite{P,RM,RG}, we have that the BSDE~\eqref{bsde1} 
gives rise to a dynamic coherent risk-measure as follows:
%
%\smallskip
%
\begin{Def}For a given coherent driver function $g$, 
	the corresponding dynamic coherent risk measure $\rho^g=(\rho^g_{t})_{t\in[0,T]}$, 
	$\mc L^2 \to \mathcal S^2$  
	is given by $$\rho^g_t(X) = Y_t,$$
	where $(Y,Z)\in\mc S^2\times\tilde{\mc H}^2$ solves \eqref{bsde1}. 
\end{Def}
\begin{Rem}\label{rem:spectal}
	\noindent{\bf (i)} Let $L^{\mathtt d}=(L^{\mathtt d}_t)_{t\in[0,T]}$ be given by $L^{\mathtt d}_t = \mathtt{d} t + L_t$ 
	for some $\mathtt{d}\in\R^{k\times 1}$.
	For random variables $X\in\mc L^2$ of the form $X=f(L^{\mathtt d}_T)$ 
	for some function $f:\R^k\to\R$ the dynamic coherent risk-measure $\rho^{g}(X)$ 
	is related to the following semi-linear PIDE (denoting $\dot{v}=\frac{\partial v}{\partial t}$):
	\begin{eqnarray}\label{pide1}
	\dot{v}(t,x) + \mathcal Gv(t,x) + g(t,\mathcal Dv_{t,x}) &=& 0, \q (t,x)\in[0,T)\times\R^k,\\
	v(T,x) &=& -f(x),\q x\in\R^k, \label{pide2}
	\end{eqnarray} 
	where $\mathcal D v_{t,x}:\R^k\to\R$ and $\mathcal Gv(t,x)$ are given by $\mathcal Dv_{t,x}(y)= v(t,x + y) - v(t,x)$ 
	and 
	\begin{equation*}
	\mathcal Gv(t,x) = {\mathtt d}^{\intercal} \nabla v(t,x) + \int_{\R^{k}\backslash\{0\}}\left[\mathcal Dv_{t,x}(y) - \nabla v(t,x)^{\intercal} y\right]\nu(\td y),
	\end{equation*}
	where $\nabla v = (\frac{\partial v}{\partial x_1}, \ldots, \frac{\partial v}{\partial x_k})^{\intercal}$.
	Specifically, if $v\in C^{1,1}([0,T]\times\R^k)$ solves \eqref{pide1}--\eqref{pide2} 
	such that $\nabla v(t,x)$ is bounded (uniformly in $(t,x)\in[0,T]\times\R^k$) then we have
	the stochastic representation
	\begin{eqnarray}\label{FK1}
	\rho^g_t(X) &=& \E\left[\left. - f(L^{\mathtt d}_T)  + \int_t^T g(t,Z_t)\td s \right|\F_t\right]\, =\, v(t,L^{\mathtt d}_t),\\
	Z_t(x) &=& v(t,L^{\mathtt d}_{t-} + x) - v(t,L^{\mathtt d}_{t-}), \q x\in\R^k, \label{FK2} 
	\end{eqnarray}
	with $L^{\mathtt d}_{0-}=L^{\mathtt d}_0$. This non-linear Feynman-Kac result is shown 
	by an application of It\^{o}'s lemma.
	
	\noindent{\bf (ii)} The risk measure $\rho^g$ admits a dual representation 
	\begin{equation}\label{delbaen}
	\rho^g_t(X) = \mathrm{ess.}\, \sup_{\mathbb Q\in\mathcal S^g}\E^{\mathbb Q}[-X|\mathcal F_t]
	\end{equation}
	for a certain representing subset $\mathcal S^g$ of the set $\mathcal M_1$ of probability measures that are absolutely continuous with respect to $\mathbb P$. 
	The set $\mathcal S^g$ is convex and closed (see \cite{D2006} and \cite[Corollary 12]{RG}).  
\end{Rem}

We describe next a representation result for dynamic risk-measure $\rho^g$ in terms of
the representing processes $(H^\xi)$ of the stochastic logarithms 
of the Radon-Nikodym derivatives $\xi\in L^2_+(\F_T)$ of the measure 
$\mathbb Q^\xi\in\mathcal M_2$, which are given by
\begin{equation}
\xi= \mathcal E\left(M^\xi\right)_T, \quad
M^\xi_\cdot =  
\int_{[0,\,\cdot\,]\times\mbb R^{k}\backslash\{0\}}H^\xi_s(x)\,\tilde N(\td s\times\td x),
\end{equation}
where $\mathcal E(\cdot)$ denotes the Dol\'{e}ans-Dade stochastic 
exponential. We call a $
\mathcal{B}(\mathbb{R}^d)\otimes \mathcal{U}$-measurable set $C=(C_t)_{t\in[0,T]}$
convex, closed or bounded if, for any $t\in[0,T]$, $C_t$ is convex, closed or bounded.  	

\begin{Thm}\label{thm:dual} Let $g$ be a coherent driver function. 
	Then for some $\mathcal{P}\otimes \mathcal{U}$-measurable set $C^g$ that is closed, convex, contains $0$ and is bounded, 
	we have 
	for any $t\in[0,T]$ that $\rho^g_t(X)$ satisfies \eqref{delbaen} with
	\begin{eqnarray}
	\mathcal S^g =\left\{\mathbb Q^\xi\in\mathcal{M}_{1}:\   H^\xi_s\in C^g_s\,\mbox{ for all }s\in[0,T]\right\}.
	\label{delbaen2}
	\end{eqnarray}
	Furthermore, the driver function $g$ satisfies for $(t,z)\in[0,T]\times\mathcal L^2(\nu)$
	\begin{equation}\label{g-rep}
	g(t,z) = \sup_{h\in C^g_t} 
	\int_{\mathbb R^k\backslash\{0\}} z(x) h(x)\nu(\td x).
	\end{equation}
\end{Thm}
The proof of Theorem~\ref{thm:dual} follows by a straightforward adaptation of 
the arguments given in \cite{D2006}, and is omitted. 
\begin{Remark}\label{rem:dual}
	\noindent{\bf (i)} Note that two driver functions $g_1$ and $g_2$ are equal if and only if 
	the corresponding sets $C^{g_1}$ and $C^{g_2}$  in the representation \eqref{g-rep} are equal.
	
	\noindent{\bf (ii)} Let  
	$\bar C$ be a $\mathcal U$-measurable  subset of $\mathcal L^2(\nu)$. 
	If $C^g_t = \bar C$ for all $t\in[0,T]$ 
	then the corresponding driver function is given by
	$g(t,z) = \bar g(z)$ where 
	\begin{equation}\label{g1g2}
	\bar g(z) = \sup_{k\in \bar C} \int_{\mathbb R^k\backslash\{0\}}  z(x) k(x)\nu(\td x),\quad 
	z\in\mathcal L^2(\nu).
	\end{equation}
\end{Remark}
\subsection{Convergence}\label{sec:conv}
We next turn to the question of the convergence of a sequence $(\rho^{g^{(\pi)}, (\pi)})_{\pi}$ of dynamic coherent risk measures as in Definition~\ref{def:disrim} when the mesh size $\Delta=\Delta_\pi$ tends to zero.   
Let us suppose that $(\rho^{g^{(\pi)}, (\pi)})_{\pi}$ are driven by the random walks $(L^{(\pi)})_\pi$ 
that are defined as follows:
\begin{equation}
\Delta L^{(\pi)}_t = J_t \sqrt{\Delta},\q J_t \stackrel{\text{IID}}{\sim}(p_j^{\Delta}, j\in\mathbb Z^k), \q t\in\pi\backslash\{T\},\label{Lpi}
\end{equation}
for some probability distribution $(p_j^{\Delta}, j\in\mathbb Z^k)$ on $\mathbb Z^k$ that is given as follows
in terms of a constant $c\ge 1$ (that will be specified shortly)
and a partition $(B^\Delta_j, j = (j_1,\ldots, j_k) \in \mbb Z^k)$ of 
$(\sqrt{\Delta}\mbb Z)^k$ into block sets of the form 
$$B^\Delta_j  = \prod_{j_i} A_{j_i}^\Delta,$$ where 
$A_k^\Delta = [k\sqrt{\Delta}, (k+1)\sqrt\Delta)$ if $k>0$, 
$A_k^\Delta = ((k-1)\sqrt{\Delta}, k\sqrt\Delta]$ if $k<0$ and
$A_0^\Delta = (-\sqrt\Delta, \sqrt\Delta)$:
\begin{eqnarray}\label{Lspec1}
p^{\Delta}_j &=& \nu(B^{\Delta}_j)\,\Delta ,\q j\in \mbb Z^k\backslash C_\Delta,\\ 
p^{\Delta}_j &=& 0,\q\q\q j\in C_\Delta\backslash\{0\},\\
p^{\Delta}_0 &=& 1- \sum_{j\neq 0}p^{\Delta}_j,\ \text{where}\ 
\label{Lspec2}\\
C_\Delta\, 
&=& \{j\in\Z^k: |j|\leq \sqrt{c_\Delta\nu_{2}}\}, 
\ \  c_\Delta = c + (\log(\Delta))^-, 
\label{Cdelta}
\end{eqnarray}
where, as before, $\nu_2 = \int_{\mathbb R^k\backslash\{0\}}|x|^2\nu(\td x)$.

When $\Delta\searrow 0$, we have by the dominated convergence theorem that
\begin{equation}\label{convL}
\E\left[\left(L_T^{(\pi),r} + L_T^{(\pi),s}\right)^2\right] 
\longrightarrow T \int_{\R^k\backslash\{0\}}[(x_r +x_s)^2]\nu(\td x),\q r,s \in \{1, \ldots, k\},
\end{equation}
where $L_T^{(\pi),m}$ and $x_m$, $m\in\{1,\ldots,k\}$, denote the $m$th 
coordinates of $L_T^{(\pi)}$ and $x\in\R^k$.

Moreover, we have by functional weak convergence theory 
(see {\em e.g.}  \cite[Theorem VII.3.7]{JS}) 
\begin{equation}\label{LpL}
L^{(\pi)} \stackrel{\mathrm{d}}{\longrightarrow} L,\q \text{as $\Delta\searrow 0$},
\end{equation}
where $\stackrel{\mathrm{d}}{\longrightarrow}$ denotes convergence in law in the Skorokhod $J_1$-topology 
on the space $\mathbb D([0,T],\mathbb R^k)$ of $\mathbb R^k$-valued RCLL functions.

On a suitably chosen 
probability space $L_T^{(\pi)}$ 
converges to $L_T$ in probability as $\Delta\searrow 0$. The latter convergence also holds in a 
stronger sense thanks to moment-conditions satisfied by $L^{(\pi)}_T$ that we show next. 
We define the value of $c$ as follows in terms of $\varepsilon_0>0$ given in Assumption~\ref{As:LevyC}:
\begin{equation}\label{c}
c = \sup_{x,y\in\R^k}\frac{|x+y|^{2 + \varepsilon_0}\vee 1}{(|x|^{2 + \varepsilon_0}\vee 1)(|y|^{2 + \varepsilon_0}\vee 1)},
\end{equation}
where $x\vee y = \max\{x,y\}$ for $x,y\in\mathbb R$.%

\begin{Lemma}\label{Lem:mom}
	The collection $(L^{(\pi)})_\pi$ of random walks defined in \eqref{Lpi} and \eqref{Lspec1}--\eqref{Lspec2}
	is such that we have, for  any uniform partition $\pi$ and $t\in\pi\backslash\{T\}$, 
	$\E\left[\left|\Delta L^{(\pi)}_t\right|\right]/\sqrt{\Delta} \to 0$ as $\Delta\searrow 0$ and
	\begin{eqnarray}\label{mmc}
	\E\left[\left|\Delta L^{(\pi)}_t\right|^{2+\varepsilon_0}\right]\leq \nu_{2+\varepsilon_0}\,\Delta, 
	\qquad \P\left(\left|\Delta L^{(\pi)}_t\right| = 0 \right) \ge 1 - \frac{1}{c_\Delta},
	\end{eqnarray}
	where $\varepsilon_0>0$ and $\nu_{2+\varepsilon_0}$ are as in Assumption~\ref{As:LevyC}, 
	and $c_\Delta$ is given in \eqref{Cdelta} and \eqref{c}. 
	Furthermore, we have 
	\begin{equation}\label{suppibound}
	\sup_{\pi: \Delta_\pi\in\R_+\backslash\{0\}}\E\left[\left|L^{(\pi)}_T\right|^{2+\varepsilon_0}\right]\in\R_+.
	\end{equation}
\end{Lemma}
\begin{Rem}
	Under the bound in the right-hand side of \eqref{mmc} we have numerical 
	stability of the solutions to sequence of BS$\Delta$Es driven by $(L^{(\pi)})$ (see \cite[Theorem 3.4]{MPS1}). 
\end{Rem}

\begin{proof}{}\ Letting $\pi=\pi_\Delta$ denote the partition with mesh $\Delta\in\R_+\backslash\{0\}$ and 
	$\varepsilon=\varepsilon_0$, a first observation is that, 
	for any $t\in\pi\backslash\{T\}$, $a\in\mbb R_+\backslash\{0\}$ and $p\in[2,2+\varepsilon]$, we have 
	by Chebyshev's inequality
	\begin{eqnarray}
	\P(|\Delta L^{(\pi)}_t| > a) 
	\label{est-c2}
	&\leq&\nu(z\in\mbb R^k: |z|\ge a)\,\Delta \\
	&\leq& \frac{\nu_{p}}{a^{p}}\,\Delta,
	\label{est-c}
	\end{eqnarray}
	where, as before, $\nu_p = \int_{\R^k\backslash\{0\}} |x|^p \nu(\td x)$. By multiplying \eqref{est-c2} 
	by $p\,a^{p-1}$ and integrating we have the estimate 
	\begin{equation}\label{key}
	\E[|\Delta L^{(\pi)}_t|^p]\leq \nu_p\, \Delta, \quad p\in[2,2+\varepsilon].
	\end{equation}
	Taking in \eqref{est-c} $p=2$ and $a=b\sqrt{c_\Delta\nu_2\Delta}$ and (a) 
	setting $b=1$ shows that 
	\begin{equation}\label{bbb1}
	\P(|\Delta L^{(\pi)}_t| > 0) = \P(|\Delta L^{(\pi)}_t|> \sqrt{c_\Delta\nu_2\Delta})\leq c_\Delta^{-1}
	\end{equation}
	which yields the bound in the right-hand side of \eqref{mmc}, 
	while (b) integrating over $b\ge 1$ 
	shows that $\E[|\Delta L^{(\pi)}_t|]/\sqrt{\Delta} \leq \sqrt{{\nu_2/c_\Delta}}$, which tends to zero as $\Delta\searrow 0$ in view of the form of $c_\Delta$. 
	
	To establish~\eqref{suppibound} the proof next proceeds analogously as that of the moment 
	result for L\'{e}vy processes (see \cite[Theorem 25.3]{Sato}).
	The key step is to transfer the uniform estimate of moments of the increments to a uniform estimate 
	of moments of the random walk at $T$ is the following estimate for a sub-multiplicative functions $g$
	(a function $g:\R^k\to\R$ is called sub-multiplicative, we recall, if for some $b_g\in\R_+$ and any $x,y\in\R^k$ 
	we have $g(x+y)\leq b_g\, g(x) g(y)$):
	\begin{eqnarray}\nonumber
	\E\left[g\left(L_T^{(\pi)}\right)\right] &=& \E\left[g\left(\sum_{t\in\pi\backslash\{T\}}\Delta L_t^{(\pi)}\right)\right]\\
	&\leq& b_g^{N-1} \E\left[g\left(\Delta L_{t_1}^{(\pi)}\right)\right]^N,
	\label{gnest}
	\end{eqnarray}
	where we used that the increments $\Delta L_t^{(\pi)}$, $t\in\pi\backslash\{T\}$, are independent. 
	For any $a\in\R_+$ the function $g_a$ given by $g_a(x) := |x|^{2+\varepsilon}\vee a$, we recall from \cite[Proposition 25.4]{Sato} 
	is sub-multiplicative. From~\eqref{key} and~\eqref{bbb1} we have that $\E\left[g_1\left(\Delta L^{(\pi)}_t\right)\right]$ 
	is bounded above by
	\begin{eqnarray}
	\E\left[g_0\left(\Delta L^{(\pi)}_t\right)\right] + \P\left(\left|\Delta L^{(\pi)}_t\right|\in(0,1]\right)
	&\leq& \nu_{2+\varepsilon}\Delta +  c_\Delta^{-1}.\label{g1est}
	\end{eqnarray}
	Combining the bounds~\eqref{gnest} and~\eqref{g1est} with the facts that $c$ defined in \ref{c} is such that $b_{g_1}=c$ and $c\leq c_\Delta$ 
	we have for all $N\in\mathbb N$
	\begin{equation}\label{gnest2}
	\E\left[g_1\left(L_T^{(\pi)}\right)\right] \leq c^{N-1}\left(\frac{1}{c} + \nu_{2+\varepsilon}\Delta \right)^N 
	= \frac{1}{c}\, \left(1 + \frac{c\,\nu_{2+\varepsilon}\,T}{N}\right)^N. 
	\end{equation}
	As the right-hand side of \eqref{gnest2} is bounded above by $c^{-1}\, \exp(c\,\nu_{2+\varepsilon}\,T)$ 
	we have \eqref{suppibound}, and the proof is complete.
\end{proof}
The moment-conditions in Lemma~\ref{Lem:mom} carry over to those of path-functionals as follows:
\begin{Cor}\label{lem:F}
	Assume that $F:\mathbb D([0,T],\mathbb R^k)\to\R$ satisfies for some $k\in\R_+$
	\begin{equation}\label{F-lip}
	|F(\omega)|\leq k \|\omega\|_\infty\ \text{for all $\omega\in\mathbb D([0,T],\mathbb R^k)$}, 
	\end{equation}
	where $\|\omega\|_\infty = \sup_{t\in[0,T]}|\omega(t)|$ for $\omega\in\mathbb D([0,T],\mathbb R^k)$.
	Then we have uniformly over partitions $\pi=\pi_\Delta$
	\begin{equation}
	\sup_{\Delta\in\R_+\backslash\{0\}}\E\left[\left|F\left(L^{(\pi)}\right)\right|^{2+\varepsilon_0}\right]\in\R_+. 
	\end{equation}
\end{Cor}
\begin{proof}{}\  
	For any partition $\pi$, an application of Doob's inequality to the centered 
	random walk $\bar L^{(\pi)}_t = L^{(\pi)}_t - t \E[L^{(\pi)}_1]$ 
	shows that 
	\begin{equation}\label{est}
	\E\left[\sup_{t\in\pi}|\bar L^{(\pi)}_t|^{2+\varepsilon_0}\right]\leq 
	\text{const}\, \E[|\bar L^{(\pi)}_T|^{2+\varepsilon_0}].
	\end{equation}
	The assertion now follows by combining the estimate \eqref{est} with 
	\eqref{F-lip}, the triangle inequality, the convexity of $x\mapsto |x|^{2+\varepsilon_0}$ and \eqref{suppibound} in Lemma~\ref{Lem:mom}. 
\end{proof}

To guarantee that the convergence of  the random walks $(L^{(\pi)})_\pi$ carries over 
to the convergence of the corresponding BS$\Delta$Es we impose the following condition on 
the sequence of coherent driver functions $(g^{(\pi)})_\pi$ and their piecewise-constant RCLL interpolations $(\tilde g^{(\pi)})_\pi$: 
\begin{Cond}\label{cond}{\bf (i)}
	The collection of functions $(g^{(\pi)})_\pi$ is uniformly 
	Lipschitz continuous with Lip\-schitz constants $K^{(\pi)}$ such that
	$\sup_{\pi} K^{(\pi)} \in\mbb R_+$.\\
	{\bf (ii)} for any continuous function  $h$
	for which $\sup_{x\in\R^k\backslash\{0\}}|h(x)|/|x|$ is bounded and any $t\in[0,T]$ we have
	\begin{equation*}\label{gpiconv}
	\lim_{\Delta\to 0} \tilde g^{(\pi)}(t,h) = g(t,h).
	\end{equation*}
\end{Cond}
The convergence result for BS$\Delta$Es (\cite[Theorem 4.1]{MPS1}) 
is phrased as follows in the current setting:
\begin{Thm}\label{thm:convB}
	Let $g$ be a coherent driver function and let $(L^{(\pi)})_\pi$ be as in \eqref{Lpi} and \eqref{Lspec1}--\eqref{Lspec2} and suppose that the sequence of coherent driver functions $(g^{(\pi)})_\pi$ satisfies Condition~\ref{cond}. 
	If $X^{(\pi)}\in \mathcal L^2(\F^{(\pi)}_T)$ and $X\in\mathcal L^2$ are such that $X^{(\pi)}\to X$ in distribution 
	and the collection $(\{X^{(\pi)}\}^2)_{\pi}$ is uniformly integrable, then we have (with $\tilde\rho^{g^{(\pi)},(\pi)}$ 
	the piecewise constant RCLL interpolation of $\rho^{g^{(\pi)},(\pi)}$): 
	\begin{equation}\label{rgprg}
	\tilde\rho^{g^{(\pi)},(\pi)}\left(X^{(\pi)}\right) \stackrel{\mathrm{d}}{\longrightarrow} \rho^g(X),\q\text{as $\Delta\searrow 0$}.
	\end{equation}
\end{Thm}
\section{Choquet-type integrals and iterated versions}\label{sec:ChE}
\subsection{Choquet-type integrals}\label{sec:ChI}
We describe next the Choquet-type integrals that feature in the definition of dynamic spectral risk-measures 
given in the next section. We refer to \cite{Den} for a treatment of the theory of non-linear integration.
The Choquet-type integrals that we consider are given in terms of measure distortions that 
we define next.
\begin{Def}\label{def:chq.mu} 
	Let $(\mathcal U,\mbb U,\mu)$ be a measure space.
	\BI\im[{\bf (i)}] $\Gamma:[0,\mu(\U)]\to[0,\infty]$ is called a {\em measure distortion} if 
	$\G$ is continuous and increasing with $\Gamma(0)=0$. If $\Gamma(1)=1$ then $\Gamma$ is called 
	a {\em probability distortion}.
	%\smallskip
	
	%\noindent
	\im[{\bf (ii)}] $\Gamma\circ\mu: \mc U\to [0,\infty]$ 
	denotes the set-function given by
	$  
	(\Gamma\circ\mu)(A) := \Gamma(\mu(A))$ for $A\in\mc U$.
	%\smallskip
	%
	%\noindent
	\EI
\end{Def}

\smallskip
On a given measure space $(\mbb U,\mathcal U)$ a set $A\in\mbb U$ with $\mu(A)>0$ is called an atom, we recall, if 
$C\subset A$ implies $\mu(C)\in\{0,\mu(A)\}$.  We assume throughout that the measure distortions and associated measure spaces are of the following type:

\begin{As}\label{Ass2}The measure $\mu$ on $(\U,\mc U)$ is sigma-finite and has no atoms, and the 
	measure distortion $\Gamma:[0,\mu(\U))\to\R_+$ is bounded and such that 
	\begin{equation}\label{KG}
	K_{\Gamma} := \int_{(0,{\mu(\U)})}\frac{\Gamma(y)}{y\sqrt{y}}\,\td y\in\mathbb R_+.
	\end{equation}
\end{As}

The Choquet-type integrals that we consider are defined as follows:

\begin{Def}\label{def:chq.mu2} Let $(\mathcal U,\mbb U,\mu)$  be a measure space 
	and let $\Gamma_+$ and $\Gamma_-$ be associated measure distortions which satisfy Assumption~\ref{Ass2}.
	\BI
	\im[{\bf (i)}] The Choquet-type integral 
	${\mathsf C}_+^{\Gamma_+\circ\mu}:\mc L_+^2(\mbb U,\mc U,\mu)\to\mbb R_+$  
	is given by
	\begin{equation*}
	{\mathsf C}_+^{\Gamma_+\circ\mu}(f) := \int_{[0,\infty)}(\Gamma_+\circ\mu)\left(f>x\right)\td x, \q f\in \mc L^2_+(\mbb U,\mc U,\mu), 
	\end{equation*}
	where $\{f>x\} = \{z\in\mbb U: f(z) >x\}$.
	\im[{\bf (ii)}] The Choquet-type integral $\mathsf C^{\Gamma_+\circ\mu,\G_-\circ\mu}: \mc L^2(\mbb U,\mc U,\mu)\to\mbb R$ 
	is given by 
	\begin{equation}\label{CH-exp-d}
	\mathsf C^{\Gamma_+\circ\mu,\G_-\circ\mu}(f) = {\mathsf C}_+^{\Gamma_+\circ\mu}(f^+) - {\mathsf C}_+^{\Gamma_-\circ\mu}(f^-),
	\end{equation} 
	where $x^+ = \max\{x,0\}$ and $x^- = \max\{-x,0\}$ for $x\in\R$. 
	\EI
\end{Def}
\begin{Remark}\label{rem:ch-finite}
	{\bf (i)} To see that ${\mathsf C}_+^{\Gamma\circ\mu}(f)\in\R_+$ for $f\in \mc L^2_+(\mu)$ 
	and $\mu$ and $\Gamma$ satisfying Assumption~\ref{Ass2} we note that 
	by Chebyshev's inequality, monotonicty of $\Gamma$ 
	and a change of variables, we have 
	\begin{eqnarray}
	{\mathsf C}_+^{\Gamma\circ\mu}(f)  = \int_0^\infty \Gamma(\mu(f>x))\td x 
	\leq \int_0^\infty\Gamma(|f|^2_{2,\mu}/x^2)\td x = 
	K_{\G} |f|_{2,\mu},\label{cb}
	\end{eqnarray}
	if $\mu(\mathbb U)=\infty$. If $\mu(\mathbb U)<\infty$  we find 
	by a similar line of reasoning that
	${\mathsf C}_+^{\Gamma\circ\mu}(f)\leq K'_{\G} \|f\|_{2,\mu}$ 
	with $K'_{\G} = K_{\G} + \G(\mu(\mathbb U))/\sqrt{\mu(\mathbb U)}$. 

	\noindent{\bf (ii)} Taking in Definition~\ref{def:chq.mu2}, $(\mathbb U,\mc U,\mu) = (\Omega,\F_T,\mathbb P)$, 
	and taking the measure distortions $\Gamma_+$ and $\G_-$ equal to a continuous probability 
	distortion $\Psi$ and the function $\WH\Psi$ given by 
	$\WH\Psi(x) = 1 - \Psi(1-x)$ for $x\in[0,1]$, it is straightforward to check that 
	$\Psi\circ\P$ is a {\em capacity} and the Choquet-type integral 
	of $X\in\mc L^2$ in \eqref{CH-exp-d} coincides with the classical Choquet expectation corresponding to 
	$\Psi\circ\P$:
	\begin{equation}\label{cl.chint}
	{\mathsf C}^{\Psi\circ\P,\WH\Psi\circ\P}(X) = 
	\int_{0}^\infty(\Psi\circ\P)(X>x)
	\td x - \int_{-\infty}^{0}(1-(\Psi\circ\P)(X > x))
	\td x.
	\end{equation}
	Moreover, as we have $\WH\Psi(x)\leq x\leq \Psi(x)$ for $x\in[0,1]$ it follows
	\begin{equation}\label{ch-ineq}
	{\mathsf C}^{\Psi\circ\P,\WH\Psi\circ\P}(X)\ge \E[X],
	\end{equation}
	and we have equality in \eqref{ch-ineq} for all $X\in\mc L^2$ if and only if
	$\Psi(x) = \WH\Psi(x)=x$ for $x\in[0,1]$.
\end{Remark}

We record next a robust representation result for Choquet-type integrals 
that plays an important role in the sequel. 
Let ${\mc M}_{p,\mu}$, $p\ge 1$, denote the set of measures $m$ on $(\mbb U,\mc U)$ 
that are absolutely continuous with respect to a given measure $\mu$ on this space 
with Radon-Nikodym derivatives such that $\frac{\td m}{\td\mu}\in\mc L^p_+(\mu)$.

\begin{Prop}
	\label{prop:md} 
	For a given concave measure distortion $\G$ and measure $\mu$ on $(\mbb U,\mc U)$  
	satisfying Assumption~\ref{Ass2}
	define
	\begin{eqnarray}\nonumber
	{\mc M}_{1,\mu}^{\Gamma} &:=& \le\{m\in\mc M^{ac}_{1,\mu}: m(A)\leq \Gamma(\mu(A))\
	\text{for all $A\in\mc U$ with $\mu(A)<\infty$} \ri\}.
	\end{eqnarray}
	Then we have that $\mathsf C_+^{\Gamma\circ\mu}:\mc L^2_+(\mu)\to \mbb R_+$ 
	is $K_\Gamma$-Lipschitz-continuous and
	\begin{equation}\label{robrep}
	\mathsf C_+^{\Gamma\circ\mu}(f) = \sup\{m(f): m\in {\mc M}_{1,\mu}^{\Gamma}\}\q
	\text{for $f\in \mc L_+^2(\mu)$}.
	\end{equation}
	In particular, $\mathsf C_+^{\Gamma\circ\mu}$ is positively homogeneous and subadditive, 
	that is, for any $\lambda\in\R_+$ and $f,g\in\mc L^2_+$
	\begin{equation}\label{phc}
	\mathsf C_+^{\Gamma\circ\mu}(\lambda f) = \lambda\, \mathsf C_+^{\Gamma\circ\mu}(f),\q
	\mathsf C_+^{\Gamma\circ\mu}(f + g) \leq 
	\mathsf C_+^{\Gamma\circ\mu}( f) + \mathsf C_+^{\Gamma\circ\mu}(g).
	\end{equation}
\end{Prop}
\begin{proof}{ of Proposition~\ref{prop:md}.} 
	The representation in~\eqref{robrep}, we recall, is known to hold true when 
	(a) $\Gamma(1)=1$ and (b) $\mu$ has unit mass and (c)  ${\mc M}_{1,\mu}^{\Gamma}$ is replaced 
	by the set of $m\in{\mc M}_{1,\mu}^{\Gamma}$ with $m(\mathbb U)=1$ 
	(see \cite{CD} and \cite[Corollary 4.80]{FollmerS}). 
	We note that, by positive homogeneity and (a) and (b), (c) is not needed for 
	the representation in \eqref{robrep} to hold true. 
	Let $\e>0$, let $\mu$ be as given and let $m\in{\mc M}_{1,\mu}^{\Gamma}$, and denote 
	by $O_\e$, $\e>0$, a collection of sets with finite non-zero $\mu$-measure and 
	such that $O_\e\nearrow\mathbb U$. Denoting
	\begin{eqnarray*}
		&&c_\e := \mu(O_\e),\q 
		\G_\e(\,\cdot\,) := \G(c_\e\,\cdot\,),\\ 
		&&m_\e(\td x):= I_{O_\e}(|x|)m(\td x),\q 
		\mu_\e(\td x) := c_\e^{-1} I_{O_\e}(|x|)\mu(\td x),
	\end{eqnarray*}
	we thus have for any $f\in\mc L^2_+(\mu)$ that 
	\begin{eqnarray}\nonumber
	\frac{1}{\G_{\e}(1)}\, {\mathsf C}_+^{\Gamma_\e\circ\mu_\e}(f) &=& 
	\sup\left\{m(f): m\in\mc M^{\Gamma_\e/\Gamma_\e(1)}_{1,\mu_\e}\right\} \\
	\label{plus} 
	&=& \sup\left\{\frac{1}{\G_\e(1)}\, m_\e(f): m\in\mc M^\Gamma_{1,\mu}\right\}.
	\end{eqnarray}
	Since, as is readily verified by an application of the monotone convergence theorem, 
	${\mathsf C}_+^{\Gamma_\e\circ\mu_\e}(f)\nearrow {\mathsf C}_+^{\Gamma\circ\mu}(f)$ 
	and $m_\e(f)\nearrow m(f)$ as $\e\downarrow 0$, and $\G_\e(1)\in\R_+\backslash\{0\}$, 
	we obtain \eqref{robrep} by taking $\e\searrow0$ in \eqref{plus}.
	
	The positive homogeneity and convexity 
	of $\mathsf C_+^{\G\circ\mu}(f)$ as stated in \eqref{phc} 
	follow as direct consequences of the robust representation in~\eqref{robrep}.
	
	Next we turn to the proof of Lipschitz continuity. 
	We observe that the robust representation~\eqref{robrep} of $\mathsf C_+^{\G\circ\mu}$ implies that for $u,v\in \mc L_+^2(\mu)$
	\begin{equation}\label{bd}
	|\mathsf C_+^{\G\circ\mu}(u) - \mathsf C_+^{\G\circ\mu}(v)|\leq 
	|\mathsf C_+^{\G\circ\mu}(v-u)|\vee |\mathsf C_+^{\G\circ\mu}(u-v)|.
	\end{equation}
	Using next a similar estimate as in \eqref{cb},
	we note that for $m\in {\mc M}_{1,\mu}^{\Gamma}$ 
	\begin{eqnarray*}
		\left|\frac{\td m}{\td\mu}\right|^2_{2,\mu} &=& 
		\int \le(\frac{\td m}{\td\mu}\ri)^2 \td \mu \,=\, 
		\int \frac{\td m}{\td\mu} \td m \,=\, \int_0^\infty m\le(\frac{\td m}{\td\mu} > x\ri)\td x \\
		&\leq& \int_0^\infty \Gamma\le(\mu\le(\frac{\td m}{\td\mu} > x\ri)\ri)\td x 
		\,\leq\, K_\Gamma \left|\frac{\td m}{\td\mu}\right|_{2,\mu},
	\end{eqnarray*}
	which implies $\sup_{m\in {\mc M}_{1,\mu}^{\Gamma}}\left|\frac{\td m}{\td\mu}\right|_{2,\mu}\leq K_\Gamma$ 
	and hence for $u\in\mc L^2_+(\mu)$ we have
	$|\mathsf C_+^{\G\circ\mu}(u)| \leq K_\G |u|_{2,\mu}$  (by~\eqref{robrep}). 
	The latter bound together with \eqref{bd} yields the stated Lipschitz-continuity. 
\end{proof}

\subsection{Conditional and iterated Choquet integrals}\label{sec:ChC}
Analogously, we define $\F_t$-conditional Choquet-type integrals as follows:

\begin{Def}\label{def:condch}For any $t\in[0,T]$ and probability distortions $\Psi$ and $\bar\Psi$
	satisfying Assumption~\ref{Ass2} (relative to the measure $\P$ restricted to $(\Omega,\F_t)$), 
	the conditional Choquet-type integral 
	${\mathsf C}^{\Psi\circ\P,\bar\Psi\circ\P}(\,\cdot\,|\F_t):\mc L^2\to \mc L^2_t$  
	is given by
	\begin{equation*}
	{\mathsf C}^{\Psi\circ\P,\bar\Psi\circ\P}(X|\mc F_t) 
	:= \int_{\R_+}\Psi\left(\P(X^+>x|\mc F_t)\right)
	\td x - \int_{\R_+}\bar\Psi\left(\P(X^- > x|\mc F_t)\right)
	\td x
	\end{equation*} 
	for $X\in\mc L^2$, where $\{X^\pm > x\} = \{\w\in \Omega: X^\pm(\w) >x\}$.
\end{Def}
\begin{Remark}{\bf (i)} Reasoning similarly as in Remark \ref{rem:ch-finite}(i) and as 
	in the proof of Lemma~\ref{prop:md},
	we have that (a)  for any $X\in\mc L^2$,  
	${\mathsf C}^{\Psi\circ\P,\bar\Psi\circ\P}(X|\mc F_t)$ is square-integrable; 
	and (b) the map ${\mathsf C}^{\Psi\circ\P,\bar\Psi\circ\P}(\,\cdot\,|\mc F_t)$ 
	is Lipschitz-continuous on $\mc L^2$ with Lipschitz-constant $K_\Psi + K_{\bar\Psi}$ 
	(which are given by the constant $K_\Gamma$ in \eqref{KG} with  $\mu(\mathbb U)=1$ and $\Gamma$ equal to 
	$\Psi$ and $\bar\Psi$, respectively). \\
	\noindent{\bf (ii)} The conditional Choquet expectation in~\eqref{CH-exp-d} 
	of $X\in\mc L^2$ with $\bar\Psi = \WH\Psi$ 
	may equivalently be expressed as weighted integral of the conditional Expected Shortfall of $X$ at different levels. Specifically, associated to any concave probability distortion $\Psi$ is a unique Borel measure $\mu$ on $[0,1]$ defined by $\mu(\{0\}) = 0$ 
	and by $\mu(\td s) = s F(\td s)$ for $s\in(0,1]$, 
	where $F$ is the locally finite positive 
	measure given in terms of the right-derivative $\Psi'_+$ of $\Psi$ by
	$F((s,1]) = \Psi'_+(s)$  (see \cite[Theorem 4.70]{FollmerS}). 
	It is straightforward to check that $\Psi$ satisfies Assumption~\ref{Ass2} if and only if 
	\begin{equation}\label{mu}
	\int_{(0,1]} \frac{1}{\sqrt{s}}\mu(\td s)\in\R_+\backslash\{0\}.
	\end{equation}
	The conditional Choquet expectation in Definition~\ref{def:condch} 
	can then be expressed in terms of the measure $\mu$ and the $\F_t$-conditional 
	Expected Shortfall, as follows:
	\begin{equation}\label{S-def}
	{\mathsf C}^{\Psi\circ\P,\WH\Psi\circ\P}(X|\mc F_t) = 
	\int_{(0,1]}\mathsf{ES}_\lambda(-X|\F_t)\mu(\td\lambda),\q X\in\mc L^2,
	\end{equation}
	where the $\F_t$-conditional Expected Shortfall $\mathsf{ES}_\lambda(X|\F_t)$ of $X\in\mc L^2$ 
	at level $\lambda\in(0,1]$ 
	is given in terms of the $\F_t$-conditional Value-at-Risk 
	$\mathsf{VaR}_{\lambda}(X|\F_t) = \inf\{z\in\R: \P(X<-z|\F_t)<\lambda\}$ at level $\lambda$ by
	\begin{eqnarray*}
		\mathsf{ES}_\lambda(X|\F_t) = \frac{1}{\lambda}\int_0^\lambda \mathsf{VaR}_{u}(X|\F_t)\td u,\q 
		\lambda\in(0,1].
	\end{eqnarray*}
	The proof of \eqref{S-def} follows by a straightforward adaptation to the conditional setting 
	of the proof for the static setting given in F\"{o}llmer and Schied (2011). \\
	\noindent{\bf (iv)} It follows from the representation in~\eqref{S-def} that
	the collection of  the conditional Choquet expectations 
	$X\mapsto {\mathsf C}^{\Psi\circ\P,\bar\Psi\circ\P}(- X|\mc F_t)$, $t\in[0,T]$, 
	$X\in\mc L^2$, is a dynamic coherent risk measure in the sense of Definition~\ref{def:rho} 
	(with $\mathcal I=[0,T]$).
\end{Remark}

One way to define  
a sequence of conditional spectral risk-measures 
that is adapted to the filtration 
$\mathbf F^{(\pi)}=(\F_t^{(\pi)})_{t\in\pi}$ is recursive 
in terms of conditional Choquet-integrals, as follows :
\begin{Def}\label{SpecIt}
	Given a concave probability distortion $\Psi$ 
	satisfying Assumption~\ref{Ass2} and a filtration $\mathbf F^{(\pi)} = (\F_t^{(\pi)})_{t\in\pi}$   
	the corresponding iterated spectral risk measure $\mathsf S=(\mathsf S_t)_{t\in\pi}$, $\mathsf S_t:\mc L^2(\F^{(\pi)}_T)\to\mc L^2(\F_t^{(\pi)})$ is defined recursively on the grid 
	$\pi=\pi_\Delta$  by 
	\begin{eqnarray}
	\mathsf S_t(X) = 
	\begin{cases}
	\mathsf C^{\Psi\circ\P,\WH\Psi\circ\P}\left(\mathsf S_{t+1}(X)\left|\F^{(\pi)}_t\right.\right), & 
	t\in\pi\backslash\{T\};\\
	-X, & t=T.
	\end{cases}
	\end{eqnarray}
\end{Def}
The class of iterated spectral risk-measures defined as such contains in particular the 
{\em Iterated Tail Conditional Expectation} proposed in \cite{HardyWirch} 
and is closely related to the {\em Dynamic Weighted V@R} that is defined in \cite{C}
for adapted processes via its robust representation. 
As already noted in the proof of Proposition~\ref{prop:md}, in the static case such a representation 
was derived in \cite{CD} for bounded random variables;
see also \cite[Theorems 4.79 and 4.94]{FollmerS} , and see \cite{CWV} for the extension 
to the set of measurable random variables (we refer to \cite{Cohenbintree} for families of dynamic risk measure 
defined via stochastic distortion probabilties in a binomial tree setting; see \cite{CohenElliotb} 
for a general theory of finite state BSDEs). 

We show next that iterated spectral risk measures are discrete-time time-consistent dynamic coherent risk measures 
and identify the driver function of the associated BS$\Delta$E.

\begin{Prop}\label{specdriver}
	The iterated spectral risk measure $\mathsf S=(\mathsf S_t)_{t\in\pi}$ given in Definition~\ref{SpecIt}
	is a discrete-time coherent risk measure $\rho^{\bar g_\Delta,\pi}$ 
	with driver function $\bar g_\Delta$ given by
	\begin{equation}\label{gth}
	\bar g_\Delta(t,h) = \frac{1}{\Delta}\left(\mathsf C^{\Psi\circ(\nu^{(\pi)}\Delta),\WH\Psi\circ(\nu^{(\pi)}\Delta)}\left(h(f)I_{\{f\neq 0\}}\right) - 
	\Delta\,\int_{\R^k\backslash\{0\}} h(x)\nu^{(\pi)}(\td x)\right),
	\end{equation}
	where $\nu^{(\pi)}$ is defined in \eqref{eq:nupi}.
\end{Prop}
\begin{proof}{}\ 
	It follows from Proposition~\ref{prop:md} that the function $\bar g_\Delta$ defined in \eqref{gth} 
	is a coherent driver function in the sense of Definition~\ref{def:cohg} with $\mathcal I = \pi$ and 
	$\mu=\nu^{(\pi)}$. Let $X\in\mc L^2(\mathcal F^{(\pi)})$ be arbitrary  and denote by $(Y^{(\pi)},Z^{(\pi)})$ 
	the solution of the BS$\Delta$E with driver function $\bar g_\Delta$.
	To show that the dynamic coherent risk measure corresponding to 
	$\bar g_\Delta$ coincides with the spectral risk measure $\mathsf S=(\mathsf S_t)_{t\in\pi}$ 
	it suffices to verify that 
	\begin{equation}\label{ogid}
	\bar g_\Delta(t,Z^{(\pi)}_t)\, \Delta = \mathsf S_t(X) - 
	\E\left[\left.\mathsf S_{t+1}(X)\right|\F_t^{(\pi)}\right].
	\end{equation}
	Letting $t\in\pi\backslash\{T\}$ and denoting $\Delta L=\Delta L^{(\pi)}_t$,
	we note from Definition~\ref{SpecIt} and \eqref{BSDD1}  
	that $\mathsf S_t(X) - 
	\E\left[\left.\mathsf S^{\phantom{(pi)}}_{t+1}(X)\right|\F_t^{(\pi)}\right]$ 
	is equal to
	\begin{eqnarray*}%\nonumber
		\lefteqn{\mathsf C^{\Psi\circ\P,\WH\Psi\circ\P}\left(\le.\mathsf S^{\phantom{(pi)}}_{t+1}(X)\right|\F_t^{(\pi)}\right)  - 
			\E\left[\left.\mathsf S^{\phantom{(pi)}}_{t+1}(X)\right|\F_t^{(\pi)}\right]}\\
		%\nonumber
		&=& \mathsf C^{\Psi\circ\P,\WH\Psi\circ\P}\left(\left.Z^{(\pi)}_t(\Delta L)I_{\{\Delta L\neq 0\}}
		\right|\F_t^{(\pi)}\right)\\ %\nonumber
		&&  - \E\left[\left.Z^{(\pi)}_t(\Delta L)I_{\{\Delta L\neq 0\}}\right| \mathcal F^{(\pi)}_t
		\right]\\
		&=& \left(\left.\mathsf C^{\Psi\circ(\nu^{(\pi)} \Delta),\WH\Psi\circ(\nu^{(\pi)} \Delta)}\left(h(f)I_{\{f\neq 0\}}\right) - 
		\Delta\,\int_{\R^k\backslash\{0\}} h(x)\nu^{(\pi)}(\td x)\right)\right|_{h=Z^{(\pi)}_t},
		%\label{star}
	\end{eqnarray*}
	where we used that, due to stationarity of the increments of $L^{(\pi)}$, 
	$\Delta L^{(\pi)}_t$ (which has law $\nu^{(\pi)}\Delta$) is independent of $t$. 
	Thus we have \eqref{ogid} and the proof is complete.
\end{proof}
\section{Dynamic spectral risk measures}\label{sec:DSR}
With the previous results in hand we move to 
the definition of dynamic spectral risk-measures in continuous time. 
Let us fix in the sequel a pair of concave measure distortions functions
$\Gamma_+$ and $\Gamma_-$ that satisfy Assumption~\ref{Ass2} and are such that $\G_-(x)\leq x$ for $x\in\mbb R_+$.
We define dynamic spectral risk measures to be those coherent spectral risk measures $\rho^g$ for which the driver 
functions $g$ are given in terms of Choquet integrals, as follows:

\begin{Def}\label{def:distdriver}
	The {\em spectral driver function} $\bar g:\mc L^2(\nu)\to\mbb R_+$
	is given by
	\begin{eqnarray*}\label{eq:spectraldriver}
		\bar{g}(u) &:=& {\mathsf C}_+^{\Gamma_+\circ\nu}(u^+)
		+ {\mathsf C}_+^{\Gamma_-\circ\nu}(u^-) 
	\end{eqnarray*}
	\noindent for  $u\in \mc L^2(\nu)$. 
\end{Def}
By Lemma~\ref{prop:md} we have that $\bar g$ is 
Lipschitz-continu\-ous, positively homogeneous and convex, 
so that $\bar g$ is a coherent driver function in the sense of Definition~\ref{def:cohg}. 
The corresponding dynamic coherent risk-measure $\rho^{\bar g}$ is the object of study 
for the remainder of the paper, which we label as follows:
\begin{Def}
	The dynamic coherent risk-measure $\rho^{\bar g}$ with spectral driver function $\bar g$ given in Definition~\ref{def:distdriver} is called the (continuous-time) {\em {dynamic spectral risk-measure} 
		corresponding to measure distortions $\Gamma_+$ and $\Gamma_-$.}
\end{Def}
We next show that dynamic spectral risk measure admit a dual representation 
of the form~\eqref{delbaen} and \eqref{delbaen2} with a representing set that is explicitly expressed 
in terms of the measure distortions $\G_+$ and $\G_-$, as follows:
%
%\medskip
\begin{Thm}\label{prop:repp}
	Let $X\in \mc L^2$, $t\in[0,T]$ and let $\bar g$ be a spectral driver function.
	The dynamic spectral risk-measure $\rho^{\bar g}$ satisfies the dual representation 
	in \eqref{delbaen}, \eqref{delbaen2} with representing set
	$
	C^{\bar g}
	$ given by
	\begin{eqnarray}
	C^{\bar g} &=& \left\{H\in \mc L^2(\nu): \begin{array}{c}
	\text{for any $A\in\mc B(\R^k\backslash\{0\})$ with $\nu(A)<\infty$}\\
	\vspace{-0.3cm}\\
	- \G_-(\nu(A)) \leq \displaystyle\int_A H\td \nu \leq \G_+(\nu(A))
	\end{array}
	\right\}, \label{C2}
	\end{eqnarray}
	where $\int_A H\td \nu = \int_A H(x)\nu(\td x)$.
\end{Thm}
\begin{Example}
	The risk of a positive or negative jump arriving with a size larger than $a$, $a\in\R_+\backslash\{0\}$, 
	as quantified by the dynamic spectral risk measure $\rho^{\bar g}$ may be explicitly expressed in terms of 
	$\nu$, $\G_+$ and $\G_-$, as we show next. For any $a\in\R_+\backslash\{0\}$, 
	let $I(a) = I_{\{\sup_{t\in[0,T]}|\Delta L_t| \leq a\}} = \{N^a_T=0\}$, 
	$N^a_T = \#\{t\in[0,T]: |\Delta L_t| > a\}$ and 
	$\bar\nu(a) = \nu(\{y: |y|>a\})$. 
	While $\E[I(a)] = \exp(- \bar\nu(a) T)$ (since $N_T$ follows a Poisson distribution with parameter $T\bar\nu(a)$), 
	the values of $I(a)$ and $-I(a)$ under $\rho^{\bar g}$ are given as follows:
	\begin{eqnarray*}
		\rho_0^{\bar g}(I(a)) &=& - \exp(-T\{\bar\nu(a)+\G_+(\bar\nu(a))\}),\\ 
		\rho_0^{\bar g}(-I(a)) &=&  \exp(-T\{\bar\nu(a)-\G_-(\bar\nu(a))\}).
	\end{eqnarray*}
	These expressions follow by deploying the dual representation in Theorem~\ref{prop:repp}   
	and Girsanov's theorem ({\em e.g.}, Theorems III.3.24 and III.5.19 in Jacod and Shiryaev (1987)): 
	we have that $\rho^{\bar g}_0(I(a))$ is equal to
	\begin{eqnarray*}
		\sup_{\Q^\xi\in\mathcal S^g} \E^{\Q^\xi}[-I(a)] &=& \sup_{\Q^\xi\in\mathcal S^g}
		\E\left[-\exp\left(-\int_0^T \int_{(a,\infty)}(1 + H^\xi_t(y))\nu(\td y)\td t\right)\right]\\
		&=& -\exp(-T\ovl\nu(a))\exp\left(- T\G_+(\ovl\nu(a))\right),
	\end{eqnarray*}
	while the expression for $\rho_0^{\bar g}(-I(a))$ follows in a similar manner.
\end{Example}
\begin{proof}{ of Theorem~\ref{prop:repp}.}\ 
	In view of Theorem~\ref{thm:dual} and Remark~\ref{rem:dual}(i)--(ii) 
	it suffices to verify that for any  $ h\in L^2(\nu)$ we have
	\begin{eqnarray}
	\sup_{k\in C^{\bar g}} \int h\,k\,\td\nu = 
	\mathsf C_+^{\G_+\circ\nu}(h^+) + \mathsf C_+^{\G_-\circ\nu}(h^-), 
	\label{2ii}
	\end{eqnarray}
	where $\int h\,k\,\td\nu = 
	\int_{\R^k\backslash\{0\}}  h(x) k(x)\nu(\td x)$.
	
	Our next observation is that the set $C^{\bar g}$ in \eqref{C2} admits the following equivalent representation:
	\begin{eqnarray}\label{C2a}
	C^{\bar g} = \left\{U\in \mc L^2(\nu): 
	\begin{array}{c}
	\text{for any $A\in\mc B(\R^k\backslash\{0\})$ with $\nu(A)<\infty$}\\
	\vspace{-0.3cm}\\
	\displaystyle\int_A U^+ \td\nu \leq \G_+(\nu(A)),\ 
	\displaystyle\int_A U^- \td\nu \leq \G_-(\nu(A))
	\end{array}
	\right\}.
	\end{eqnarray}
	To see that this is the case, we note that, for any $U\in L^2(\nu)$, 
	we have $-U^-\leq U\leq U^+$, while $U^+ = U_1$ and $-U^- = U_2$
	for $U_1 = U I_{\{U\ge 0\}}$ and $U_2 = U I_{\{U< 0\}}$.
	
	To see that \eqref{2ii} holds 
	we note from \eqref{C2a}, Proposition~\ref{prop:md} and the identity
	$$ h\, k^* =  h^+\,k_1^+ +  h^-\,k_2^-,\q
	k^*= k_1^+\,I_{\{ h>0\}} -  k_2^-\,I_{\{ h<0\}},$$
	for any $ h, k_1, 
	k_2\in\mc L^2(\nu)$,
	that $\bar g(h) = \sup_{k\in C^{\bar g}} \int h\,k\,\td\nu$ is bounded below by 
	\begin{equation*}%\label{ggeq}
	\sup_{k_1,k_2\in C^{\bar g}} \int h\,k^*\,\td\nu = 
	\sup_{ k\in C^{\bar g}}\int  h^+\, k^+\,\td\nu + 
	\sup_{ k\in C^{\bar g}}\int  h^-\, k^-\,\td\nu,
	\end{equation*}
	which is by Proposition~\ref{prop:md} equal to
	$\mathsf C_+^{\G_+\circ\nu}( h^+) + \mathsf C_+^{\G_-\circ\nu}( h^-)$. 
	Given this lower bound and the fact  that $\bar g( h)$ is bounded above by
	\begin{equation*}
	\sup_{ k\in C^{\bar g}}\int  h^+\, k\,\td\nu + 
	\sup_{ k\in C^{\bar g}}\int  h^-\,(- k)\,\td\nu 
	\leq \sup_{ k\in C^{\bar g}}\int  h^+\, k^+\,\td\nu + 
	\sup_{ k\in C^{\bar g}}\int  h^-\, k^-\,\td\nu, 
	\end{equation*}
	we conclude that \eqref{2ii} holds true. 
\end{proof}

\section{Limit theorem}\label{sec:lim}
We next turn to the functional limit theorem which shows that dynamic spectral risk measures arise as a limit of iterated spectral risk measures, under a suitable scaling of the corresponding probability distortions. 
We suppose that, uniformly in $p\in[0,1]$, $\Psi_\Delta(p)-p$ 
scales in the mesh size $\Delta$ and the measure distortions $\G_+$ and $\G_-$ as follows:
$$
\Psi_\Delta(p) = p + \Delta\left\{\Gamma_+(p/\Delta)I_{[0,\frac{1}{2}]}(p) 
+ \Gamma_-((1-p)/\Delta)I_{(\frac{1}{2},1]}(p)\right\} + \mathrm{o}(\Delta)\ (\Delta\searrow 0).
$$
Specifically, the condition that we require is phrased as follows: 
\begin{Def}\label{def:PsiD} 
	We denote by $(\Psi_\Delta)_{\Delta\in(0,1]}$ a sequence of probability distortions 
	that is such that $\Psi_\Delta$ and $\WH\Psi_\Delta$ given by $\WH\Psi_\Delta(p) = 1 - \Psi_\Delta(1-p)$ satisfy
	Assumption~\ref{Ass2} with respect to the measure $\mu(\td x)\equiv \P(\Delta L^{(\pi)}_{t_1}\in\td x)$ 
	and we have
	\begin{eqnarray}\label{limit}
	\lim_{\Delta\searrow 0} \Upsilon_\Delta = 0,\q \Upsilon_\Delta = \sup_{x\in(0,1)} \left|\frac{\Psi_\Delta(x) - x}{\G_\Delta(x)\Delta} - 1\right|,
	\end{eqnarray}
	where for $\Delta\in(0,1]$ and $x\in[0,1]$ 
	$$
	\G_\Delta(x) = \Gamma_+(x/\Delta)I_{[0,\frac{1}{2}]}(x) + \Gamma_-((1-x)/\Delta)I_{(\frac{1}{2},1]}(x). 
	$$
	Here, we recall, $\G_+$ and $\G_-$ denote the given concave measure distortions 
	which are such that $\G_-(x)\leq x$ for $x\in\R_+$ and  
	Assumption~\ref{Ass2} holds with $\mu(\td x)\equiv\nu(\td x)$ and $\Gamma\equiv \G_+$ or $\G_-$. 
\end{Def}
The functional limit result is phrased as follows in terms of 
the sequence of piecewise-constant RCLL extensions $(\tilde L^{(\pi)})_{\pi}$ of 
the random walks $(L^{(\pi)})_{\pi}$ 
given by 
$$\tilde L^{(\pi)}_t := L^{(\pi)}_{\Delta^{-1}[t\Delta]}, \quad t\in[0,T],
$$ 
where $[r]=\sup\{n\in\mathbb N\cup\{0\}: n\leq r\}$ for $r\in\R_+$.
\begin{Thm}\label{main}
	Given a sequence of probability distortions $(\Psi_\Delta)_{\Delta\in(0,1]}$ as in Definition~\ref{def:PsiD} 
	and given filtrations $\mathbf F^{(\pi)} = (\F_t^{(\pi)})_{t\in\pi}$,
	let $\mathsf S^{\Delta}=(\mathsf S^{\Delta}_t)_{t\in\pi}$, $\Delta\in(0,1]$, denote the corresponding iterated 
	spectral risk-measures as given in Definition~\ref{SpecIt} and let $\bar g$ denote 
	the spectral driver function from Definition~\ref{def:distdriver}. 
	Let the set of $\omega\in\mathbb D([0,T],\mathbb R^k) $ at which  $F:\mathbb D([0,T],\mathbb R^k)\to\R$ is discontinuous in the Skorokhod $J_1$-topology 
	be a null-set under the law of $L$  
	and assume that for some $k\in\R_+$
	\begin{equation}\label{Flip}
	|F(\omega)|\leq k \|\omega\|_\infty\ \text{for all $\omega\in\mathbb D([0,T],\mathbb R^k)$}, 
	\end{equation}
	where $\|\omega\|_\infty = \sup_{t\in[0,T]}|\omega(t)|$ for $\omega\in\mathbb D([0,T],\mathbb R^k)$.
	Then we have  
	\begin{equation}\label{speco}
	\tilde{\mathsf S}^{\Delta}\left(F\left(\tilde L^{(\pi)}\right)\right) \stackrel{\mathrm{d}}{\longrightarrow} 
	\rho^{\bar g}\left(F\left(L\right)\right),\quad \Delta\searrow 0, 
	\end{equation}
	where $\tilde{\mathsf S}^{\Delta}_t = {\mathsf S}^{\Delta}_{\Delta^{-1}[t\Delta]}$, $t\in[0,T]$.
\end{Thm}
\begin{Rem}\label{remlim}
	\noindent{\bf (i)} Given two concave probability distortions $\Psi_+$ and $\Psi_-$ 
	satisfying the integrability condition~\eqref{KG} (with $\mu(\mathbb U)=1$) one may explicitly construct 
	a sequence $(\Psi_\Delta)_{\Delta\in(0,1]}$ satisfying Definition~\ref{def:PsiD} as follows:
	\begin{eqnarray*}%\label{PSD}
		\Psi_\Delta(p) =  p + (\G_+(p/\Delta)I_{[0,\frac{1}{2}]}(p) + \G_-((1-p)/\Delta)I_{(\frac{1}{2},1]}(p))\,\Delta,\q p\in[0,1],
	\end{eqnarray*}
	where, inspired by \cite{Eb}, we suppose that 
	the functions $\G_+, \G_-:\R_+\to\R_+$ are given by 
	\begin{eqnarray*}%\label{G+G-}
		\G_+(x) = a\,\Psi_+(1 - \te{-c x}),\q \G_-(x) = \frac{b}{d}\,\Psi_-(1 - \te{-d x}),\q x\in\R_+,
	\end{eqnarray*}
	for some $a$, $b$, $c$ and $d\in\R_+\backslash\{0\}$  satisfying 
	the restrictions 
	\begin{equation}\label{abcd}
	\G_+(1/(2\D)) = \G_-(1/(2\D)) < 1/(2\Delta),\quad b\,\Psi'_-(0^+)\in(0,1),
	\end{equation}
	where $f'(0^+)$ denote the right-derivative of a function $f$ at $x=0$.
	It is straightforward to check that, for any $\D\in(0,1]$, 
	$\Psi_\D$ is a concave probability distortion (the first condition in \eqref{abcd} guarantees 
	continuity at $p=1/2$ and $\Psi_\Delta(1/2)<1$) and that 
	$\G_-(x)\leq x$  for any $x\in\R_+$ (as consequence of the second condition in \eqref{abcd}).
	Furthermore, we have that the limit in \eqref{limit} holds. 
	
	\noindent{\bf (ii)} 
	Examples of functionals $F$ that satisfy condition~\eqref{Flip} 
	include (a) a European call option payoff with strike $K\in\R_+$ ($F(\w) = (\w(T)-K)^+$); (b) 
	the time-average ($F(\w) = \frac{1}{T}\int_0^T \w(s)\td s$) and 
	(c) the running maximum ($F(\w) = \sup_{s\in[0,T]}\w(s)$).
	
	\noindent{\bf (iii)} We note that $\Upsilon_\Delta$ may be equivalently expressed in terms 
	of $\Psi_\Delta$ and $\WH\Psi_\Delta$ as follows:
	$$
	\Upsilon_\Delta = \sup_{x\in(0,\frac{1}{2}]}\left|\frac{\Psi_\Delta(x) - x}{\G_+(x/\Delta)\Delta} - 1\right|
	\,\bigvee\,\sup_{x\in(0,\frac{1}{2})}\left|\frac{x - \WH\Psi_\Delta(x)}{\G_-(x/\Delta)\Delta} - 1\right|.
	$$
	
	\noindent{\bf (iv)} We next provide an example to show the necessicity of scaling the probability distortions.
	For a given uniform partition $\pi=\pi_\Delta$ of $[0,T]$ with mesh $\Delta$, a probability distortion $\Psi$ 
	and $a_+,a_-\in\R_+\backslash\{0\}$, let us consider the risk-charge under the iterated spectral risk measure 
	$\mathsf S$ corresponding to $\Psi$ of the following statistic $X^{(\pi)}$ of the jump-sizes of 
	$L^{(\pi)} = (L^{(\pi),1},\ldots, L^{(\pi),k})$:
	\begin{equation}\label{nanb}
	X^{(\pi)} := N^+_\pi - N^-_\pi, \q
	N^\pm_\pi = \#\left\{t\in\pi\backslash\{T\}: \sum_{i=1}^k |\Delta L^{(\pi),i}_{t}|^\pm > a_\pm\right\}. 
	\end{equation}
	From the form~\eqref{Ypi}--\eqref{tildeZpi} of the solution of the BS$\Delta$E associated to 
	the iterated spectral risk measure $\mathsf S$ we have that $Z^{(\pi)}$ is given by
	\begin{eqnarray}
	&&Z^{(\pi)}_t(x) = z^{(\pi)}_+(x) - z^{(\pi)}_-(x), \q z^{(\pi)}_\pm(x) = I_{A_\pm}(x),\\
	&&A_\pm = \left\{z\in\R^k\backslash\{0\}: \sum_{i=1}^k|z_i|^\pm > a_\pm\right\}.
	\end{eqnarray}
	As a consequence, we have from \eqref{gth} in Proposition~\ref{specdriver} that
	the driver function takes the form
	\begin{eqnarray*}
		\bar g_{\Delta}(t,Z^{(\pi)}_t)\Delta &=& 
		\mathsf C^{\Psi\circ(\nu^{(\pi)}\Delta),\WH\Psi\circ(\nu^{(\pi)}\Delta)}\left(z^{(\pi)}_+(f)-z^{(\pi)}_-(f)\right) \\ && - 
		\Delta\,\int_{\R^k\backslash\{0\}} (z^{(\pi)}_+(x)-z^{(\pi)}_-(x))\nu^{(\pi)}(\td x)\\
		&=& \Psi(\P(\Delta L^{(\pi)}_{t_1}\in A_+)) - \P\left(\Delta L^{(\pi)}_{t_1}\in A_+\right)\\
		&&  +\, \P\left(\Delta L^{(\pi)}_{t_1}\in A_-\right) - \WH\Psi(\P(\Delta L^{(\pi)}_{t_1}\in A_-)).
	\end{eqnarray*}
	For given $t\in\pi\backslash\{T\}$ 
	the iterated spectral risk-measure $\mathsf S_t(X^{(\pi)})$, 
	may thus be expressed as follows in terms of the functions $D_\D^+$ and $D_\D^-:[0,\Delta^{-1}]\to\R_+$ given by
	$D_\Delta^+(x) = \Psi(x\,\Delta) - x$ and $D_\Delta^-(x) = x - \WH\Psi(x\,\Delta )$:
	\begin{eqnarray*}
		\lefteqn{\mathsf S_t(X^{(\pi)}) - \E[X^{(\pi)}|\F^{(\pi)}_t] = \E\left[\left.\sum_{s\ge t, s\in\pi\backslash\{T\}}\bar g_{\Delta}(s,Z^{(\pi)}_s)\Delta \right|\F^{(\pi)}_t \right]}\\
		&=& (T-t)\left( \frac{1}{\Delta}\,D^+_\Delta(\Delta^{-1}\,\P(\Delta L^{(\pi)}_{t_1}\in A_+)) 
		+ \frac{1}{\Delta}\,D^-_\Delta(\Delta^{-1}\,\P(\Delta L^{(\pi)}_{t_1}\in A_-))\right).
	\end{eqnarray*}
	Note that, as $\Delta\searrow 0$, 
	$\Delta^{-1}\,\P(\Delta L^{(\pi)}_{t_1}\in A_\pm)\to \nu(A_\pm)$ 
	and 
	$$\E[X^{(\pi)}|\F^{(\pi)}_t]\to (T-t)(\nu(A_+)-\nu(A_-)) + N_t^+ - N_t^-,$$
	where $N^\pm_t = \#\{s\in(0,t]: L_s-L_{s-}\in A_\pm\}$. Hence, 
	this suggests that for 
	the sequence of iterated spectral risk-measures to converge, 
	$\Delta^{-1}\,D_\Delta^+(x)$ and $\Delta^{-1}\,D_\Delta^-(x)$ are to admit limits as $\Delta\searrow 0$.
\end{Rem}

\begin{proof}{ of Theorem~\ref{main}.}
	We note first that, as $L^{(\pi)}\stackrel{\mathrm d}{\to} L$ when $\Delta\searrow 0$, 
	$F(L^{(\pi)})$ converges in distribution to $F(L)$, which is element of $\mc L^2$. 
	Furthermore, by Corollary~\ref{lem:F}, the collection 
	$\{F(L^{(\pi)})^2\}_\pi$ is uniformly integrable.	
	Thus, in view of Theorem~\ref{thm:convB} it suffices next to verify that the sequence of 
	driver functions $(\bar g_\Delta)_{\Delta\in(0,1]}$
	of the iterated spectral risk measures $\mathsf S^{\Delta}$ given in Proposition~\ref{specdriver} 
	satisfies Condition~\ref{cond}, which we proceed to do.
	
	Let $t\in[0,T]$. Our first observation is by subadditivity and nonnegativity of $\bar g_{\D}$ we have 
	for any $h,k\in L^2(\nu^{(\pi)})$ 
	\begin{equation}\label{gd}
	\left|\bar g_{\Delta}(t,h) - \bar g_{\Delta}(t,k)\right|
	\leq \bar g_{\Delta}(t,h - k) \vee \bar g_{\Delta}(t,k-h),
	\end{equation}
	so that to verify Condition~\ref{cond}(i) it suffices to show that $\bar g_{\Delta}(t,h)/|h|_{2,\pi}$ 
	is uniformly bounded. We have for any $\Delta\in(0,1]$ and $h\in L^2(\nu^{(\pi)})$ that  
	\begin{eqnarray}\nonumber
	\bar g_{\Delta}(t,h) &=& \frac{1}{\Delta}\left(\mathsf C^{\Psi_\Delta\circ(\nu^{(\pi)}\Delta),\WH\Psi_\Delta\circ(\nu^{(\pi)}\Delta)}(h) -  \Delta\int h\td \nu^{(\pi)}\right)\\ \nonumber
	&=& \frac{1}{\Delta}\left(\mathsf C_+^{\Psi_\Delta\circ(\nu^{(\pi)}\Delta)}(h^+) -  \Delta\int h^+\td \nu^{(\pi)}\right) \\ \nonumber
	&& + 
	\frac{1}{\Delta}\left(\Delta\int h^-\td \nu^{(\pi)} - \mathsf C_+^{\WH\Psi_\Delta\circ(\nu^{(\pi)}\Delta)}(h^-)\right)\\ 
	&=& \mathsf C_+^{\G_+\circ\nu^{(\pi)}}(h^+) + R^\Delta(h^+) + 
	\mathsf C_+^{\G_-\circ\nu^{(\pi)}}(h^-) + \WH R^\Delta(h^-),\label{C+C-}
	\end{eqnarray}
	where the remainder terms $R^\Delta(h^+)$ and $\WH R^\Delta(h^-)$ 
	are given as follows in terms of the identity function $I:[0,1]\to[0,1]$, $I(x)=x$:
	\begin{eqnarray*}
		R^\Delta(h^+) &=& \frac{1}{\Delta}\int_0^\infty \left[(\Psi_\Delta - I)\left(\nu^{(\pi)}(h^+ > x)\Delta\right) 
		- \G_+\left(\nu^{(\pi)}(h^+ > x)\right)\Delta\right]\td x\\
		\WH R^\Delta(h^-) &=& \frac{1}{\Delta}\int_0^\infty \left[(I-\WH\Psi_\Delta)\left(\nu^{(\pi)}(h^- > x)\Delta\right) 
		- \G_-\left(\nu^{(\pi)}(h^- > x)\right)\Delta\right]\td x.
	\end{eqnarray*}
	Since by Chebyshev's inequality $\nu^{(\pi)}(h^\pm>x)\leq |h^\pm|^2_{2,\nu^{(\pi)}}/x^2$ for 
	$x\in\R_+\backslash\{0\}$, it follows that for $$x\ge H^\pm := |h^\pm|_{2,\nu^{(\pi)}}\sqrt{2\Delta}$$
	the mass of $\Delta\, \nu^{(\pi)}(h^\pm>x)$ is bounded above by $1/2$. Recalling the form of 
	$\Upsilon_\Delta$ (see Remark~\ref{remlim}(iii)) and that $\G_++ \G_-$ is bounded (by $\G_\infty$ say)
	we have
	\begin{eqnarray}\nonumber
	|R^\Delta(h^+)| &\leq& \Upsilon_\Delta\, 
	\mathsf C_+^{\G_+\circ\nu^{(\pi)}}(h^+)
	\\ && \nonumber + \int_0^{H^{+}}(\G_+(\nu^{(\pi)}(h^+>x)) + \G_-(\nu^{(\pi)}(h^-\leq x)) \td x\\ \label{Rd1}
	&\leq& \Upsilon_\Delta\, \mathsf C_+^{\G_+\circ\nu^{(\pi)}}(h^+) 
	+  \, H^+ \, \G_\infty,\\
	|\WH R^\Delta(h^-)| &\leq& \Upsilon_\Delta\, \mathsf C_+^{\G_-\circ\nu^{(\pi)}}(h^-)
	+  H^- \,\G_\infty.
	\label{Rd2}
	\end{eqnarray}
	Combining \eqref{gd}, \eqref{C+C-}, \eqref{Rd1} and \eqref{Rd2} and the $K_{\G_+}$-
	and $K_{\G_-}$-Lipschitz-continuity of $\mathsf C_+^{\G_+\circ\nu^{(\pi)}}$ and $\mathsf C_-^{\G_-\circ\nu^{(\pi)}}$ (Proposition~\ref{prop:md}) and the fact that the values $\mathsf C_+^{\G_+\circ\nu^{(\pi)}}(0)$ and $\mathsf C_+^{\G_-\circ\nu^{(\pi)}}(0)$ are equal to $0$, 
	we find 
	\begin{equation}\label{gup}
	\left|\bar g_{\Delta}(h)\right| \leq  \tilde C\,|h|_{2,\nu^{(\pi)}},
	\end{equation}
	where $\tilde C = (K_{\G_+} + K_{\G_-} + 2\sqrt{2}\, \G_\infty)(1 + \sup_{\Delta\in(0,1]}\Upsilon_\D)$ 
	is finite by the limit \eqref{limit} 
	in Definition~\ref{def:PsiD}. This completes the proof of Condition~\ref{cond}(i).
	
	We turn next to the proof of Condition~\ref{cond}(ii). 
	Let $h$ be a continuous function that is such that $c_h := \sup |h(x)/x|\in\mathbb R_+$.
	Since $\nu^{(\pi)}$ converges weakly to $\nu$, we have 
	that $\nu^{(\pi)}(h > x) \to \nu(h > x)$ at  $x\in\R_+\backslash\{0\}$ that are points of continuity.
	Hence, as $\G_\pm$ are continuous it follows that $\G_\pm(\nu^{(\pi)}(h > x)) \to \G_\pm(\nu(h > x))$ at such $x$. 
	Next we show that the latter functions are dominated by an integrable function. 
	By Chebyshev's inequality, $\G_\pm(\nu^{(\pi)}(h > x)) \leq \G_\pm(|h|^2_{2,\nu^{(\pi)}}/x^2)$ 
	while it follows from the inequality \eqref{key} that $\nu_2^{(\pi)}\leq \nu_2$, 
	where $\nu_2^{(\pi)} = \int_{\mathbb R^k\backslash\{0\}} |x|^2\nu^{(\pi)}(\td x)$.
	Hence we have the bound
	\begin{eqnarray*}
		|h|_{2,\nu^{(\pi)}} \leq c_h\, \sqrt{\nu_2^{(\pi)}} 
		\leq c_h\, \sqrt{\nu_2}.
	\end{eqnarray*}
	Also, for any $d\in\mathbb R_+$,  $\G_\pm(d^2/x^2)$ is integrable:
	\begin{eqnarray*}
		\int_0^\infty \G_\pm(d^2/x^2)\td x = K_{\G_\pm}\, d,
	\end{eqnarray*}
	where $K_{\Gamma_\pm}$ is given in \eqref{KG}. 
	As a consequence, from the dominated convergence theorem we have that
	$\mathsf C_+^{\G_\pm\circ\nu^{(\pi)}}(h^\pm)\to \mathsf C_+^{\G_\pm\circ\nu}(h^\pm)$
	as $\Delta\searrow 0$. 
	Further, in view of \eqref{limit}, $R^\Delta(h^+)$ and $\WH R^\Delta(h^-)$ tend to zero as $\Delta\searrow 0$.
	This establishes Condition~\ref{cond}(ii), and the proof is complete.
\end{proof}
\section{Dynamically optimal portfolio allocation}\label{sec:DynOpt}
We next consider dynamic portfolio problems concerning balancing gain and 
risk as quantified by the DSR. We suppose the investment horizon is equal to $T>0$ and 
consider the DSR associated to the spectral driver function $\bar g$.
In this section we impose the following restriction on the L\'{e}vy measure $\nu$:
\begin{As}
	The support of $\nu$ is included in the set $(-1,\infty)^k$.
\end{As}
We suppose that the financial market consists of a risk-free bond  
and $n$ risky stocks with discounted prices $\hat S = (\hat S^1, \ldots, \hat S^n)$ 
evolving according to the following system of SDEs:
\begin{eqnarray*}
	\frac{\td\hat S^i_t}{\hat S^i_{t-}} &=& \mathtt{d}^i\td t + \int_{\R^{k}\backslash\{0\}} {\mathtt R}^i\, x\, 
	\tilde N(\td t\times\td x), 
	\q i=1,\ldots, n, t\in(0,T],\\ 
	\hat S_0 &=& s_0\in(\R_+\backslash\{0\})^k,
\end{eqnarray*}
where $\mathtt{d}^i\in\R$ is the excess log-return and ${\mathtt R}^i\in\R^{k}$ 
is the (row) vector of jump-coefficients with non-negative coordinates that are
such that $(({\mathtt R}^i)^{\intercal} \mbf 1 \leq 1$ (where $\mathbf 1\in\R^{k}$ denotes the $k$- column vector of ones and where, for any vector $v$, $v^\intercal$ denotes its transposition). 
Given the form of the model we have $\hat S^i_t\in\mc L^2_t$ 
and $\hat S^i_t>0$ for any $i=1,\ldots, k$ and $t\in[0,T]$.

Let us consider the case of a small investor whose trades have a negligible impact on the price and 
let us adopt the classical frictionless and self-financing 
setting (no transaction cost, infinitely divisible assets, continuous-time trading, 
no funds are infused into or withdrawn from the portfolio at intermediate times, {\em etc.}). 
At any time $t\in[0,T]$ the investor decides 
to allocate part $\theta_t^i$ of the current wealth 
for investment into the stock $\hat S^i$, $i=1,\ldots, n$, so that, if $X^{\theta}_{t-}$ 
denotes the discounted wealth just before time $t$, we have that $\theta_t^i X^{\theta}_{t-}/\hat S^i_{t-}$ 
is the number of stocks $i$ held in the portfolio at time $t$. 
We suppose that certain limits are placed on the leverage ratio of the portfolio and 
on the size of the short-holdings in the various stocks, and that this restriction 
is phrased in terms of a bounded and closed set $\mathcal B\subset\R^n$ as the requirement that
\begin{equation}\label{thetaB}
\theta_t(\omega)\in\mathcal B\ \ \text{for any $(t,\omega)\in[0,T]\times\Omega$}.
\end{equation}
\vspace{-0.6cm}
\begin{Example}\label{ex:B}
	To impose constraints on the fractions of the current wealth 
	invested in the bond account and the stock accounts we take
	$$\mathcal B = \left\{x\in(\R_+)^n: x_i\ge -L_i, \sum_{i=1}^n x_i\leq 1+L_0\right\}
	$$ for some $L_0, \ldots, L_n\in\R_+$. In particular, by taking $L_i>0$ we impose a limit on the
	borrowing ($i=0$)
	or the number of stock $i$ that may be shorted ($i\neq 0$). 
	The case of a ``long only'' investor
	that has no short-sales and only invests own wealth (no borrowing) corresponds to 
	taking in $L_0=L_1=\cdots=L_n=0$.
\end{Example}

We call an allocation strategy $\theta=(\theta_t)_{t\in[0,T]}$ {\em admissible} if 
$\theta$ is predictable and \eqref{thetaB} holds. 
We denote by $\mathcal A$ the collection of admissible allocation strategies.
Denoting by ${\mathtt R}=({\mathtt R}^i)_{i=1,\ldots, k}$ the $\R^{n\times k}$-matrix 
with $i$th row equal to ${\mathtt R}^i$, we have that
the discounted value $X^\theta=(X^\theta_t)_{t\in[0,T]}$ of a portfolio 
corresponding to $\theta\in\mathcal A$ evolves according to the SDE
\begin{eqnarray*}%\label{SDEX1}
	\frac{\td X^{\theta}_t}{X^{\theta}_{t-}} &=& \theta_t^{\intercal}\, {\mathtt d}\td t +   
	\int_{\R^k\backslash\{0\}} \theta_t^{\intercal}{\mathtt R}\, x\, \tilde N(\td t \times\td x),\q t\in(0,\tau^\theta\wedge T],\\
	X^{\theta}_0 &=& x \in \R_+\backslash\{0\},\q X^\theta_t = X^\theta_{\tau^\theta\wedge T}, \ t\in(\tau^\theta\wedge T,T],
	%\label{SDEX2}
\end{eqnarray*}
where $\tau^{\theta} = \inf\{t\in[0,T]: X^{\theta}_t < 0\}$ (with $\inf\emptyset=+\infty$)
is the first time that the value of the portfolio 
becomes negative, when the investor has to stop trading.

\subsection{Portfolio optimisation under dynamic spectral risk measures}\label{cDSR}
We consider next the stochastic optimisation problem given in terms of DSR by the following criterion that is to be minimised for $t\in[0,T]$:
\begin{equation}\label{minp}
\tilde {\mathcal J}_t^\theta = \rho^{\bar g}_t(X^{\theta}_{T\wedge\tau^\theta}), 
\end{equation}
The investor's problem is to identify a 
stochastic process $\tilde{\mathcal J}^*=(\tilde{\mathcal J}^*_t)_{t\in[0,T]}$
and an allocation strategy $\theta^*\in{\mathcal A}$ such that
\begin{equation}\label{Jstar}
\tilde{\mathcal J}^*_t = \text{ess.\,}\inf_{\theta\in{\mathcal A}} \tilde{\mathcal J}_t^\theta 
= \tilde{\mathcal J}^{\theta^*}_t,\q t\in[0,T]. 
\end{equation}
While the problem in~\eqref{Jstar} 
may be solved via a BSDE approach 
(as used in for instance \cite{B,R} to analyse utility optimisation 
and robust portfolio choice prolems), due to its Markovian nature 
it may also be approached via classical methods based on an associated 
Hamilton-Jacobi-Bellman equation---this is the method that we present here. One class of allocation strategies are those of feedback-type  that are defined as follows.

\begin{Def}\label{fbeq}
	Denote by $\tilde\Theta$ the set of functions $\bar\theta:[0,T]\times\R_+\to\mathcal B$ that are
	such that the following SDE admits a unique 
	solution $X^{\bar\theta} = (X^{\bar\theta}_t)_{t\in[0,T]}$:
	\begin{eqnarray}\label{barX1}
	\frac{\td X^{\bar\theta}_t}{X^{\bar\theta}_{t-}} &=& \bar\theta(t,X^{\bar\theta}_{t-})^{\intercal}\mathtt d\,\td t + 
	\bar\theta(t,X^{\bar\theta}_{t-})^{\intercal}\, {\mathtt R}\, x\,\tilde N(\td t\times\td x),\q t\in(0,\tau^{\bar\theta}],\\
	X^{\bar\theta}_0 &=& x,\qquad X^{\bar\theta}_t=X^{\bar\theta}_{\tau^{\bar\theta}\wedge T},\ t\in(\tau^{\bar\theta}\wedge T, T],\label{barX2}
	\end{eqnarray}
	where $\tau^{\bar\theta}=\inf\{t\in[0,T]: X^{\bar\theta}_t < 0\}$.
	A strategy $\theta\in\mathcal A$ is called a {\em feedback allocation strategy} if 
	there exists a feedback function $\bar\theta\in\bar\Theta$ such that  
	$$\theta_t = \bar\theta(\tau^{\bar\theta}\wedge t,X^{\bar\theta}_{\tau^{\bar\theta}\wedge (t-)}), \quad t\in[0,T],
	$$ 
	where $X^{\bar\theta}_{0-}=X^{\bar\theta}_{0}$ and 
	$X^{\bar\theta}$ solves the SDE in \eqref{barX1}--\eqref{barX2}. 
\end{Def}

Associated to a given allocation strategy of feedback-type $\bar\theta\in\bar\Theta$ there exists a 
value function $V^{\bar\theta}$ satisfying $J^*_t = V^{\bar\theta}(t,X^{\bar\theta}_t)$ 
for $t\in[0,T]$ (as a consequence of the Markov property). 
If sufficiently regular, the function 
$V^{\bar\theta}$ satisfies a semi-linear PIDE that is given in terms of certain operators $\mathcal D^{\theta}$ and ${\mathcal G}^\theta$ indexed by $\theta\in\mathcal B$. 
For any function $f\in C^{1,1}([0,T]\times\R)$, these operators
are equal to the functions $\mathcal D^{\theta}_{t,x}f:\R^k\to\R$ and
$\mathcal G^\theta f:[0,T]\times\R_+\backslash\{0\}\to\R$ 
that are given in terms of
\begin{equation}
{\mathtt d}_{\theta} = \theta^\intercal \mathtt d,\qquad
{\mathtt R}_\theta = \theta^{\intercal}{\mathtt R}, \qquad \theta\in\mathcal B,
\label{dt-rt}
\end{equation}
by (denoting $f' = \frac{\partial f}{\partial x}$)
\begin{eqnarray*}
	(\mathcal D^{\theta}_{t,x}f)(y) &=& f(t,x+x\,{\mathtt R}_\theta\, y) - f(t,x)\\
	{\mathcal G}^\theta f(t,x) &=& {\mathtt d}_{\theta} f'(t,x) + \int_{\mathbb R^{k\times 1}\backslash\{0\}} 
	\le\{(\mathcal D^{\theta}_{t,x}f)(y) - f'(t,x)\,x\, {\mathtt R}_\theta\, y \ri\}\nu(\td y).
\end{eqnarray*}
The non-linear Feynman-Kac formula 
(see Remark~\ref{rem:spectal}) implies that if the following semi-linear PIDE 
has a sufficienly regular solution it is equal to $V^{\bar\theta}$:
\begin{eqnarray*}
	&&\dot{v}(t,x) + \mathcal G^{\bar\theta(t,x)} v(t,x) +  \bar g\left(\mathcal D^{\bar\theta(t,x)}_{t,x} v\right) = 0,\quad 
	(t,x)\in[0,T)\times\R_+\backslash\{0\},\\ 
	&&v(t,x) = - x, \q (t,x)\in[0,T)\times(\R\backslash\R_+)\cup\{0\},\\
	&&v(T,x) = - x,\q x\in\R.
\end{eqnarray*}
Standard arguments suggest then that 
if the optimal allocation strategy $\theta^*$ is of feedback-type and the corresponding 
value-function $V$ is sufficiently regular, then $V$ satisfies the following Hamilton-Jacobi-Bellman (HJB) 
equation:
\begin{eqnarray}\label{HJB1b}
&&\!\!\!\!\!\!\!\!\!\!\!\! \dot{V}(t,x) + \inf_{\theta\in\mathcal B}\left\{\mathcal G^{\theta}V(t,x) + 
\bar g(\mathcal D^{\theta}_{t,x} V)\right\} = 0, 
\quad (t,x)\in[0,T)\times \R_+\backslash\{0\},\\
\label{HJB3b}
&&\!\!\!\!\!\!\!\!\!\!\!\! V(t,x) = - x, \qquad\qquad t\in[0,T)\times (\R\backslash\R_+)\cup\{0\},\\
&&\!\!\!\!\!\!\!\!\!\!\!\! V(T,x) = - x, \qquad\qquad x\in\R.
\label{HJB4b}
\end{eqnarray}
Next we verify that a sufficiently smooth solution of 
the HJB equation gives rise to a solution of the optimisation problem in~\eqref{Jstar}. Let $C_b^{1,1}([0,T]\times\R)$ denote the set of $C^{1,1}$-functions $f:[0,T]\times\R\to\R$ with 
bounded first-order derivatives.
\begin{Thm}\label{thm:verif}
	Let $w\in C^{1,1}_b([0,T]\times\R)$ 
	be a solution of the HJB-equation ~\eqref{HJB1b}--\eqref{HJB4b} 
	and let the function $\tilde\theta:[0,T]\times\R_+\to\mathcal B$, $(t,x)\mapsto\tilde\theta(t,x)$ given by
	$$\tilde\theta(t,x) \in \mathrm{arg.}\,\sup_{\theta\in\mathcal B}\left[{\mathcal G}^{\theta}w(t,x) + 
	\bar g(\mathcal D^{\theta}_{t,x} w)\right]
	$$
	be such that $\tilde\theta\in\bar\Theta$.
	Then the feedback strategy $\tilde\theta^* = (\tilde\theta^*_t)_{t\in[0,T]}$ 
	with feedback function $\tilde\theta$ is optimal for \eqref{Jstar} and we have
	$
	\tilde{\mathcal J}^*_t = \tilde{\mathcal J}_t^{\tilde\theta^*} = 
	w(t, X^{\tilde\theta}_{t\wedge\tau^{\tilde\theta}}),
	$
	where $X^{\tilde\theta}$ solves the SDE in \eqref{barX1} and \eqref{barX2} 
	with $\bar\theta$ replaced by $\tilde\theta$. 
\end{Thm}
\begin{proof}{}\,
	Letting $\theta\in{\mathcal A}$ be an arbitrary admissible strategy, $t < \tau^\theta\wedge T$ and
	$w$ as stated in the theorem, 
	we find by an application of It\^{o}'s lemma that
	\begin{eqnarray}\nonumber
	\lefteqn{w\left(T\wedge\tau^{\theta}, X^{\theta}_{T\wedge\tau^{\theta}}\right) - 
		w(t,X^{\theta}_t) + \int_t^{T\wedge\tau^{\theta}}\bar g(\mathcal D^{\theta_s}w_{s,X^{\theta}_s})\td s}\\
	&&\!\!\!\!\!=
	\int_t^{T\wedge\tau^{\theta}}\left\{
	\dot{w} + {\mathcal G}^{\theta_s}w\right\}\left(s,X^{\theta}_s\right)
	+ \bar g(\mathcal D^{\theta_s}w_{s,X^{\theta}_s})\td s
	+ M^{\theta}_{T\wedge\tau^{\theta}} - M^{\theta}_t,
	\label{ito}
	\end{eqnarray}
	where $M^{\theta}$ is the square-integrable martingale given by
	\begin{eqnarray*}
		\lefteqn{M^{\theta}_t = \int_0^t w'(s,X^{\theta}_{s-})\left(\td X^{\theta}_s - {\mathtt d}_{\theta_s} X^{\theta}_s\td s\right)}\\
		&& + \int_0^t\int_{\R^k\backslash\{0\}}\left(\mathcal D^{\theta_s} w_{s,X^{\theta}_s}(y) - w'(s,X^{\theta}_s)\, x\,
		{\mathtt R}_{\theta_s}\, y\right)\tilde N(\td s\times\td y).
	\end{eqnarray*}
	Note that by the HJB equation~\eqref{HJB1b} the first term on the right-hand side of \eqref{ito} is non-positive.
	Hence by taking conditional expectation in \eqref{ito} and using \eqref{HJB3b}--\eqref{HJB4b} 
	we have that  
	\begin{equation}\label{wineq}
	w(t,X^{\theta}_t) \leq 
	\E\left[\left.-X^{\theta}_{T\wedge\tau^{\theta}} + \int_t^{T\wedge\tau^\theta}\bar g\left(\mathcal D^{\theta_s}w_{s,X^{\theta}_s}
	\right)\td s\right|\F_t\right] = \mathcal J^{\theta}_t.
	\end{equation}
	Since $\theta\in\mathcal A$ is arbitrary we have that 
	\begin{equation}\label{wineq2}
	w(t,X^{\bar\theta}_t)\leq \mathrm{ess.}\,\inf_{\theta\in\mathcal A} \mathcal J^{\theta}_t = \mathcal J^*_t.
	\end{equation} 
	If we choose $\theta=\tilde\theta^*$, we note that the first term on the right-hand side of 
	\eqref{ito} vanishes and the inequalities in \eqref{wineq}--\eqref{wineq2} become equalities, so that 
	$\mathcal J^*_t = w(t,X^{\tilde\theta^*}_t)$. As the process $X^{\tilde\theta^*}$ coincides with the process 
	$X^{\tilde\theta}$ solving the SDE in~\eqref{barX1}---\eqref{barX2}, the proof is complete.
\end{proof}
\subsubsection{Case of a ``long-only'' investor}
We next restrict to the case of the ``long-only'' investor
(see Example~\ref{ex:B}). In this case we note that 
for any admissible allocation strategy $\theta\in\mathcal A$ the solvency constraint 
$X^{\theta}_t \in \R_+$ is satisfied for all $t\in[0,T]$ so that $\tau^{\theta}=\infty$ a.s.
We identify the optimal strategy as follows: 
\begin{Thm}\label{thm:main-l}
	Let $\theta^*\in\mathcal B$ satisfy
	\begin{equation}\label{thst}
	\theta^* \in \mathrm{arg.}\,\sup_{\theta\in\mathcal B}\{ {\mathtt d}_\theta - \bar g(-{\mathtt R}_\theta I)\},
	\end{equation}
	where ${\mathtt d}_\theta$ and ${\mathtt R}_\theta$ are given in \eqref{dt-rt} and 
	$I:\R^{k}\to\R^{k}$ is given by the column vector $I(y)=y$. Then $\tilde\theta^*=(\tilde\theta^*_t)_{t\in[0,T]}$ 
	given by $\tilde\theta^*_t\equiv\theta^*$ is an optimal strategy 
	and 
	\begin{equation}\label{JJ}
	\mathcal J^*_t = - X_t^{\theta^*}\,\exp\left((T-t)\left\{{\mathtt d}_{\theta^*} - \bar g(-{\mathtt R}_{\theta^*}\, I)\right\}\right).
	\end{equation}
\end{Thm}
\begin{proof}{}\,
	The assertions follow by an application of the verification theorem (Theorem~\ref{thm:verif}). 
	
	We note first that as the function $\theta\mapsto{\mathtt d}_\theta - \bar g(-{\mathtt R}_\theta\, I)$ 
	is concave it attains its maximum on the compact set $\mathcal B$. Thus, the set in \eqref{thst} is not empty 
	and $\theta^*$ is well-defined. Moreover, given the positive homogeneity of $g$
	it is straightforward to verify that the function $C:[0,T]\to\R$ given by
	$$C(t) = -\exp\left((T-t)\left\{{\mathtt d}_{\theta^*} - \bar g(-{\mathtt R}_{\theta^*} I)\right\}\right)$$
	satisfies the ODE 
	\begin{eqnarray*}
		&&\dot C(t) + \inf_{\theta\in\mathcal B}\left\{{\mathtt d}_{\theta}C(t) + \bar g(C(t)\,{\mathtt R}_{\theta}\, I)\right\} = 0,
		\q t\in[0,T),
		\\
		&&C(T) = -1.
	\end{eqnarray*}
	As a consequence, we have that the candidate value function $V:[0,T]\times\R_+\to\R_+$ given by 
	$V(t,x) = C(t) x$ 
	satisfies the HJB equation \eqref{HJB1b}---\eqref{HJB4b} (here we used again the positive homogeneity of $\bar g$).
	The assertions follow now by an application of Theorem~\ref{thm:verif}.
\end{proof}

\section*{Acknowledgements}
	We thank two anonymous referees 
	and an Associate Editor for careful reading and useful comments. We thank 
	H. Albrecher, J. Blanchet, P. Jevti\'{c}, R. Laeven, 
	H. Schumacher, J. Sekine and 
	participants of the London-Paris Bachelier Workshop (Paris), 
	Workshop on Mathematical Finance and Related Issues (Osaka), SF@W seminar (Warwick), 
	Workshop on Advanced Modelling in Mathematical Finance (Kiel), 
	and Workshop on Models and Numerics in Financial Mathematics (Leiden) for useful discussions and suggestions. 
	A previous version of the paper was entitled 
	``On consistent valuations based on distorted expectations: 
	from multinomial walks to L\'{e}vy processes''.


\begin{thebibliography}{100}
	\bibitem{Ac} Acerbi, C.  Spectral measures of risk: A coherent representation 
	of subjective risk aversion. {\it Journal of Banking and Finance} 26, 1505--1518, 2002.
	\bibitem{AT} Acerbi, C. and Tasche, D.  On the coherence of expected shortfall. 
	{\it Journal of Banking and Finance} 26, 1487–-1503, 2002.
	\bibitem{ADEH} Artzner P., Delbaen F., Eber J.M. and Heath D.
	Coherent measures of risk. {\it Mathematical Finance} 9,  203--228, 1999.
	\bibitem{ADEHK} Artzner P., Delbaen F., Eber J.M., Heath D. and Ku, H.
	Coherent multiperiod risk adjusted values and Bellman’s
	principle. {\it Annals of Operations Research} 152, 5–-22, 2007.
	\bibitem{BBP} Barles, G., Buckdahn, R. and Pardoux, E. 
	Backward stochastic differential equations and
	integral-partial differential equations. {\it Stochastics and
		Stochastics Reports} 60,  57--83, 1997.
	\bibitem{B} Becherer, D.  Bounded solutions to Backward SDE's with jumps 
	for utility optimization and indifference hedging. 
	{\em Annals of Applied Probability} 16, 2027--2054, 2006.
	\bibitem{BN2} Bion-Nadal, J. Dynamic risk measures: Time consistency 
	and risk measures from BMO martingales. {\em Finance and Stochastics} 12,  219--244, 2008.
	\bibitem{BN1} Bion-Nadal, J.  Time consistent dynamic risk measures 
	{\em Stochastic Processes and Their Applications} 119,  633--654, 2009.
	\bibitem{CD} Carlier, G., Dana, R. A.  Core of convex
	distortions of a probability. {\it J. Economic Theory} 113,  199--222, 2003. 
	\bibitem{CE} Chen, Z.  
	and Epstein, L.  Ambiguity, risk, and asset returns in continuous 
	time. {\em Econometrica} 70, 1403--1443, 2002.
	\bibitem{CDK} Cheridito, P., Delbaen, F. and Kupper, M.  Dynamic monetary risk
	measures for bounded discrete time processes.
	{\em Elect. J. Probab.} 11, 57--106, 2006.    
	\bibitem{CWV} Cherny, A.S. Weighted V$@$R and its properties. 
	{\it Finance and Stochastics} 10,  367--393, 2006.    
	\bibitem{C} Cherny, A.S.  Capital Allocation  and Risk Contribution
	with Discrete-Time Coherent Risk. {\em Mathematical Finance} 19, 13--40, 2009.
	\bibitem{CohenElliotb} Cohen, S. and Elliott, R.  A
	general theory of finite state backward stochastic
	difference equations. {\it Stochastic Processes and Their Applications} 
	120, 442--466, 2010.
	\bibitem{cohenelliota} Cohen, S. and Elliott, R. 
	Backward Stochastic Difference Equations and Nearly Time-Consistent Nonlinear Expectations.
	{\em SIAM Journal of Control and Optimization} 49, 125-–139, 2011.
	\bibitem{CHMP} Coquet, F., Hu, Y., Memin, J. and Peng, S.
	Filtration-consistent nonlinear expectations and
	related g-expectations. {\it Prob. Th. Rel. Fields} 123,  1--27, 2002.
	\bibitem{CZ} Czichowsky, C. Time consistent mean-variance portfolio selection in 
	discrete and continuous time. 
	{\em Finance and Stochastics} 17, 227--271, 2013.    
	\bibitem{D2006} Delbaen, F. The Structure of m–Stable Sets and in Particular of the Set of Risk Neutral Measures.
	{\em Lecture Notes in Mathematics} 1874, 215--258, 2006. 
	\bibitem{Den} Denneberg, D.  {\em Non-additive measure
		and integral}, Kluwer Academic Publishers, 1994.
	\bibitem{DNM} Dolinsky, Y., Nutz, M. and Soner, H. M. 
	Weak Approximation of $G$-expectation. {\em Stochastic
		Processes and Their Applications} 122,  664--675, 2012.
	\bibitem{DE1992}
	Duffie, D. and Epstein, L.~G.  Stochastic differential utility. 
	{\em Econometrica} 60, 353--394, 1992.
	\bibitem{Eb} 
	Eberlein, E., Madan, D.B., Pistorius, M. and Yor, M.  
	Bid and Ask Prices as Non-Linear Continuous Time G-Expectations 
	Based on Distortions. {\em Mathematics and Financial Economics} 8, 265--289, 2014.
	\bibitem{Cohenbintree} 
	Elliott, R.J., Siu, T.K. and Cohen, S.N. 
	Backward Stochastic Difference Equations for Dynamic Convex Risk 
	Measures on a Binomial Tree. {\it Journal of Applied Probability} 52(3), 2015.
	\bibitem{ES2003} Epstein, L.G. and Schneider, M. 
	Recursive multiple-priors. {\em Journal of Economic Theory} 113, 1-–31, 2003.
	\bibitem{EZ1989}
	Epstein, 
	L.G. and Zin, S.E. Substitution, risk aversion, and the temporal 
	behavior of consumption and asset returns: a theoretical framework. 
	{\em Econometrica} 57, 937-–969, 1989.
	\bibitem{FollmerP}F\"{o}llmer, H. and Penner, I. 
	Convex risk measures and the dynamics of their penalty functions.
	{\em Statistics and Decisioncs} 24, 61--96, 2006. 
	\bibitem{FollmerS} F\"{o}llmer, H. and Schied, A.  {\it
		Stochastic Finance. } De Gruyter, 3rd edition, 2011.
	\bibitem{HS2008} Hansen, L. P. and 
	and Sargent, T.J.  {\em Robustness}. Princeton University Press, 2008.
	\bibitem{HardyWirch} Hardy, M.R. and Wirch, J.L.  
	The Iterated CTE: a Dynamic Risk Measure. 
	{\it North American Actuarial Journal} 8, 62--75, 2004.  
	\bibitem{JS} Jacod, J. and Shiryaev, A.N.  {\it Limit
		Theorems for Stochastic Processes}, Springer, 1987.
	\bibitem{JR} Jobert, A. and Rogers, L.  Valuations and
	dynamic convex risk measures. {\it Mathematical Finance}
	18,  1--22, 2008.
	\bibitem{KlSchw} Kl\"{o}ppel, S. and Schweizer, M. 
	Dynamic Indifference Valuation via Convex Risk Measures.
	{\em Mathematical Finance} 17, 599--627, 2007.
	\bibitem{K1960} Koopmans, T.C.  Stationary ordinal 
	utility and impatience. {\em Econometrica} 28, 287–-309, 1960.
	\bibitem{KP1978}
	Kreps, M.K. and Porteus, E.L.  Temporal resolution of uncertainty and dynamic 
	choice theory. {\em Econometrica} 46, 185-–200, 1978.
	\bibitem{Kusuoka} Kusuoka, S.  On law-invariant coherent risk measures.
	{\it Advances in Mathematical Economics} 3,  83--95, 2001.
	\bibitem{KuMo} Kusuoka, S. and Morimoto, Y.  
	Homogeneous Law-Invariant Coherent Multi-period Value Measures and their Limits.
	{\it J. Math. Sc.  Univ. Tokyo} 14, 117--156, 2007.
	\bibitem{R} Laeven, R.J.A. and Stadje, M.  
	Robust Portfolio Choice and Indifference Valuation. 
	{\em Mathematics of Operations Research} 39, 1109--1141, 2014.
	\bibitem{MPS1} Madan, D., Pistorius, M. and Stadje, M. 
	Convergence of BS$\Delta$Es driven by random walks to BSDEs: the case of (in)finite activity jumps with general driver.
	{\em Stochastic Processes and Their Applications} 126, 1553--1584, 2016.% {\tt doi:10.1016/j.spa.2015.11.013}
	\bibitem{P} Peng, S.  A generalized dynamic programming
	principle and Hamilton-Jacobi-Bellman equation {\it
		Stochastics and Stochastics Reports} {38},  119--134, 1992.
	\bibitem{PengR1} Peng, S. Nonlinear Expectations, Nonlinear Evaluations, 
	and Risk Measures. Stochastic Methods in Finance. In: Fritelli, M., Runggaldier, W.
	(Eds.) Lecture Notes in Mathematics, Springer, 165--254, 2004. 
	\bibitem{Riedel} Riedel, F.  Dynamic coherent risk measures.
	{\it Stochastic Processes and Their Applications} {112}, 185--200, 2004.
	\bibitem{RoorS}Roorda, B. and  Schumacher, J.M. 
	Time consistency conditions for acceptability measures - with an application to Tail Value at Risk. 
	{\em Insurance Mathematics and Economics} 40, 209--230, 2007. 
	\bibitem{RG} Rosazza Gianin, E.  Risk measures via $g$-expectations.
	{\it Insurance Mathematics and Economics} 39,  19--34, 2006. 
	\bibitem{RM} Royer, M.  Backward stochastic differential
	equations with jumps and related non-linear expectations.
	{\it Stochastic Processes and Their Applications} 116, 
	1358--1376, 2006.
	\bibitem{Sato} Sato, K.  {\it L\'{e}vy processes and
		infinitely divisible distributions}, Cambridge University
	Press, Cambridge, 1999.
	\bibitem{SIME} Stadje, M.  Extending dynamic convex risk measures from 
	discrete time to continuous time: A convergence approach
	{\it Insurance Mathematics and Economics} 47,  391--404, 2010.
	\bibitem{Strotz} Strotz, R.H. Myopia and inconsistency in Dynamic Utility Optimisation.
	{\em Review of Economic Studies} 23, 165--180, 1955. 
	\bibitem{Tutsch} Tutsch, S.  Update rules for convex risk measures.
	{\em Quantitative Finance} 8, 833--843, 2008.
	\bibitem{Wang96} Wang, S.  Premium calculation by transforming the layer premium
	density, {\it ASTIN Bulletin} 26,  71--92, 1996.
	\bibitem{Weber} Weber, S.   Distribution-Invariant Risk Measures, Information, and Dynamic Consistency. 
	{\em Mathematical Finance} 16, 419--442, 2006.
\end{thebibliography}


\end{document}